\documentclass{article}
\usepackage{amssymb}
\usepackage{amsfonts}
\usepackage{amsmath}

\newtheorem{theorem}{Theorem}

\newtheorem{corollary}[theorem]{Corollary}

\newtheorem{lemma}[theorem]{Lemma}

\newtheorem{remark}[theorem]{Remark}

\begin{document}

\title{Random walks conditioned to stay nonnegative and branching processes
in nonfavorable random environment}
\author{Congzao Dong\thanks{%
Xidian University, 266 Xinglong Section of Xifeng Road, Xi'an, Shaanxi,
710126, China Email: czdong@xidian.edu.cn}, Elena Dyakonova\thanks{%
Steklov Mathematical Institute of Russian Academy of Sciences, 8 Gubkina
St., Moscow 119991 Russia Email: elena@mi-ras.ru}, and Vladimir Vatutin
\thanks{%
Steklov Mathematical Institute of Russian Academy of Sciences, 8 Gubkina
St., Moscow 119991 Russia Email: vatutin@mi-ras.ru}}
\maketitle

\begin{abstract}
Let $\left\{ S_{n},n\geq 0\right\} $ be a random walk whose increments
belong without centering to the domain of attraction of an $\alpha $-stable
law $\{Y_{t},t\geq 0\}$, i.e. $S_{nt}/a_{n}\Rightarrow Y_{t},t\geq 0,$ for
some scaling constants $a_{n}$. Assuming that $S_{0}=o(a_{n})$ and $%
S_{n}\leq \varphi (n)=o(a_{n}),$ we prove several conditional limit theorems
for \ the distribution of $S_{n-m}$ given $m=o(n)$ and $\min_{0\leq k\leq
n}S_{k}\geq 0$. These theorems complement the statements established by F.
Caravenna and L. Chaumont in 2013. The obtained results are applied for
studying the population size of a critical branching process evolving in
nonfavorable environment.
\end{abstract}

\section{Introduction\label{s1}}

The study of random walks conditioned to stay positive or nonnegative has a
long history (see, for instance, \cite{BD94},\cite{Bolt76},\cite{BrD2006},%
\cite{Chau97},\cite{Car2005},\cite{CC2008},\cite{CC2013},\cite{CD2010}, \cite%
{Don12}\cite{Ig74},\cite{Ka76}, \cite{Lig68}, \cite{VW09} to mention a few).
The reason to establish different invariance principles for conditioned
processes, to which the majority of the papers is devoted, comes not only
from the pure theoretical interest but is explained by applications in the
theory of branching processes in deterministic and random environment \cite%
{Afan2011}, \cite{KV2017}, \cite{VD2017}, in statistical physics, with
particular reference to random polymers \cite{Hol74}, and other fields.

The present paper also deals with random walks conditioned to stay positive.
Our study is motivated by Caravenna and Chaumont paper \cite{CC2013}. The
authors of \cite{CC2013} show that if a random walk has increments with
distribution belonging (without centering) to the domain of attraction of a
stable law, then its excursion of length $n$ conditioned to stay positive
and suitably rescaled, converges in distribution toward the excursion of the
corresponding stable Levy process conditioned to stay positive. We
complement their result by considering the behavior of the properly scaled
trajectory of the random walk in a left vicinity $n-m\,\ $of the excursion,
where $m=o(n)$. We show that, depending on the rate of decay the ratio $m/n$
to zero, three different limiting distributions appear. These regimes
results in three different growth rates of a critical branching process
evolving in an nonfavorable random environment and surviving up to a distant
moment $n$. Our statements for branching processes in random environment
complement the respective theorems established in \cite{VD2017},\cite{VD2022}
and \cite{VD2023}, where, given the process does not extinct before the
moment $n$, the distribution of the population sizes were investigated for $%
m=o(n)$ and $m=\left[ nt\right] ,0<t\leq 1.$

The paper is organized as follows.

In section \ref{s2} we introduce basic notation, describe our main
assumptions for random walks, recall some known local limit theorems and
prove auxiliary statements for random walks conditioned to stay positive.

In section \ref{s3} we present a number of conditional limit theorems for the
excursions of a lattice random walk conditioned to stay positive.

In section \ref{s4} we establish the respective conditional limit theorem
for the excursions of random walks with increments having absolutely
continuous distributions.

Section \ref{s5} we prove a conditional limit theorem for sequences of
random variables convergent almost surely.

In section \ref{s6} we apply the results for the excursion of random walks
conditioned to stay positive to get a description of the population size of
a critical branching process, evolving in unfavorable environment given its
survival for a long time.

\section{Notation and assumptions\label{s2}}

We mainly use the notation and assumption introduced in \cite{CC2013} to
simplify references to the results of the paper.

In the sequel we denote by $C_{1},C_{2},...$ some absolute constants which
may not be the same in different formulas.

We set $\mathbb{N}\mathbf{:=}\left\{ 1,2,...\right\} $ and $\mathbb{N}_{0}:=%
\mathbb{N}\mathbf{\cup }\left\{ 0\right\} $. Given two positive sequences $%
\left\{ c_{n},n\in \mathbb{N}\right\} $, $\left\{ d_{n},n\in \mathbb{N}%
\right\} ,$ we write as usual $c_{n}\sim d_{n}$ if $lim_{n\rightarrow \infty
}c_{n}/d_{n}=1,c_{n}=o(d_{n})$ if $lim_{n\rightarrow \infty }c_{n}/d_{n}=0$
and $c_{n}=O(d_{n})$ if $\lim \sup_{n\rightarrow \infty }c_{n}/d_{n}<\infty $%
.

We recall that a positive sequence $\left\{ c_{n},n\in \mathbb{N}\right\} $
--- or a real function $c(x)$ --- is said to be regularly varying at
infinity with index $\gamma \in \mathbb{R}$ , denoted $\left\{ c_{n},n\in
\mathbb{N}\right\} \in R_{\gamma }$ or $c(x)\in R_{\gamma }$ if $c_{n}\sim
n^{\gamma }l(n)$ $(c(x)\sim x^{\gamma }l(x))$, where $l(x)$ is a slowly
varying function, i.e. a positive real function with the property that $%
l(cx)/l(x)\rightarrow 1$ as $x\rightarrow \infty $ for all fixed $c>0$.

Let $(\alpha ,\beta )$ be a pair of nonnegative numbers belonging to the
union of the sets
\begin{equation*}
\{0<\alpha <1;\,|\beta |<1\}\cup \{1<\alpha <2;|\beta |\leq 1\}\cup \{\alpha
=1,\beta =0\}\cup \{\alpha =2,\beta =0\}.
\end{equation*}%
For a random variable $X,$ we write $X\in \mathcal{D}\left( \alpha ,\beta
\right) $ if the distribution of $X$ belongs (without centering) to the
domain of attraction of a stable law with density $g(x)=g_{\alpha ,\beta
}(x),x\in (-\infty ,+\infty ),$ and the characteristic function%
\begin{equation*}
G\mathbb{(}w\mathbb{)}=\int_{-\infty }^{+\infty }e^{iwx}g(x)\,dx=\exp
\left\{ -c|w|^{\,\alpha }\left( 1-i\beta \frac{w}{|w|}\tan \frac{\pi \alpha
}{2}\right) \right\} ,\ c>0.
\end{equation*}%
Consider a random walk%
\begin{equation*}
S_{0}=0,\quad S_{n}=X_{1}+...+X_{n},\,n\geq 1,
\end{equation*}%
with independent and identically distributed increments. Throughout this
section we assume that the random walk $(\mathcal{S}=\left\{ S_{n},n\in
\mathbb{N}_{0}\right\} ,\mathbb{P})$ meets the following conditions.

\textbf{Hypothesis} \textbf{A1.} \textit{The increments }$X_{n},n=1,2,...$%
\textit{\ of }$\mathcal{S}$\textit{\ belong to }$D\left( \alpha ,\beta
\right) $\textit{\ with }$\left\vert \beta \right\vert <1$\textit{. }

This means, in particular, that there is an increasing sequence of positive
numbers
\begin{equation}
a_{n}\ :=\ n^{1/\alpha }\ell (n),n=1,2,...,  \label{defA}
\end{equation}%
where $\ell (1),\ell (2),\ldots $ is a slowly varying sequence, such that,
as $n\rightarrow \infty $
\begin{equation}
\mathcal{L}\left\{ \frac{S_{\left[ nt\right] }}{a_{n}},t\geq 0\right\}
\Longrightarrow \mathcal{L}\left\{ Y_{t},t\geq 0\right\} ,  \label{LimitLaw}
\end{equation}%
where the symbol $\Longrightarrow $ stands for the convergence in
distribution in the space $D[0,+\infty )$ with Skorokhod topology and the
process $(\mathcal{Y}=\left\{ Y_{t},t\geq 0\right\} ,\mathbf{P})$ is
strictly stable with index $\alpha \in $ $(0,2]$ and positivity parameter $%
\rho :=\mathbf{P}\left( Y_{1}>0\right) \in $ $(0,1)$.

\bigskip

We denote by $\Omega ^{RW}:=\mathbb{R}^{\mathbb{N}_{0}}$ the discrete path
space and by $\Omega :=D\left( [0,\infty ),\mathbb{R}\right) $ the space of
real-valued c\`{a}dl\'{a}g paths on $[0,\infty )$, equipped with the
Skorokhod topology, which turns it into a Polish space, and with the
corresponding Borel $\sigma $-field. We also set $\Omega _{N}^{RW}:=\mathbb{R%
}^{\left\{ 0,1,...,N\right\} }$ and $\Omega _{t}:=D\left( [0,t],\mathbb{R}%
\right) $ for $t\in (0,\infty ).$ We denote $\mathbb{P}_{x}$ the law of the
random walk started at $x\in \mathbb{R}$, i.e. the law on $\Omega ^{RW}$ of $%
\mathcal{S}+x$ under $\mathbb{P}$, and, for a process $\mathcal{Y}=\left\{
Y_{t},t\geq 0\right\} $ on $\Omega $ with a law $\mathbf{P}$ on $\Omega $ we
denote by $\mathbf{P}_{a}$ the law of $\mathcal{Y}+a$ under $\mathbf{P,}$
for $a\in \mathbb{R}$.

\bigskip

\textbf{Hypothesis A2.} \textit{We assume that either of the following
assumptions hold:}

\textit{\textbullet\ (}$(h;c)$\textit{-lattice case) The law of }$X_{1}$%
\textit{\ under }$\mathbb{P}$\textit{\ is supported by the lattice }$c+hZ$%
\textit{, where the span }$h>0$\textit{\ is chosen to be maximal (i.e., the
law of }$X_{1}$\textit{\ is not supported by }$c_{0}+h_{0}Z$\textit{, for
any }$h_{0}>h$\textit{\ and }$c_{0}\in R$)\textit{. Note that we may take }$%
c\in \lbrack 0;h)$\textit{.}

\textit{\textbullet\ (absolutely continuous case) The law of }$X_{1}$\textit{%
\ under }$\mathbb{P}$\textit{\ is absolutely continuous with respect to the
Lebesgue measure on }$R$\textit{, and there exists }$n\in N$\textit{\ such
that the density }$f_{n}(x):=\mathbb{P}(S_{n}\in dx)/dx$\textit{\ of }$S_{n}$%
\textit{\ is essentially bounded (therefore, }$f_{n}(x)\in L^{\infty }$%
\textit{).}

Set%
\begin{equation*}
L_{N}=\min_{1\leq k\leq N}S_{k},\quad L_{N}^{\ast }=\min_{0\leq k\leq
N}S_{k}.
\end{equation*}

Given $N\in \mathbb{N}$ and $x,y\in \lbrack 0,\infty )$, by the (law of the)
random walk bridge of length, conditioned to stay positive, starting at $x$
and ending at $y$, we mean either of the following laws on $\Omega
_{N}^{RW}: $

\begin{equation}
\mathbb{P}_{x,y}^{\uparrow ,N}\left( \cdot \right) :=\mathbb{P}_{x}\left(
\cdot \,|L_{N}\geq 0,S_{N}=y\right) .  \label{MesNonegat}
\end{equation}

Clearly, for the law $\mathbb{P}_{x,y}^{\uparrow ,N}\left( \cdot \right) $
to be well-defined in the lattice case, it is necessary to assume that the
conditioning events have positive probabilities

\begin{equation}
q_{N}^{+}\left( x,y\right) :=\mathbb{P}_{x}\left( L_{N}\geq 0,S_{N}=y\right)
>0.  \label{Posit_0}
\end{equation}

Analogously, in the absolutely continuous case we need that $%
f_{N}^{+}(x,y)>0 $, where

\begin{eqnarray}
&&f_{N}^{+}(x,y):=\frac{\mathbb{P}_{x}\left( L_{N-1}>0,S_{N}\in dy\right) }{%
dy}  \notag \\
&=&\int_{\mathcal{K}(N-1)}\left[ f(s_{1}-x)\left(
\prod_{i=2}^{N-1}f(s_{i}-s_{i-1})\right) f(s_{N-1}-y)\right]
ds_{1}...ds_{N-1}  \label{ContDistrib}
\end{eqnarray}%
with $\mathcal{K}(N-1):=\left\{ s_{1}>0,...,s_{N-1}>0\right\} $ and $f(\cdot
)=f_{1}(\cdot )$ the density of the random walk step $X_{1}$.

For $t\in \lbrack 0,\infty )$ and $a,b\in \lbrack 0,\infty )$ we denote by $%
\mathbf{P}_{a,b}^{\uparrow ,t}$ the law on $\Omega _{t}$ corresponding to
the bridge of the Levy process of length $t$ conditioned to stay positive,
starting at$\ a$ and ending at $b$, which informally may be written as
\begin{equation*}
\mathbf{P}_{a,b}^{\uparrow ,t}\left( \cdot \right) =\mathbf{P}_{a}\left(
\cdot |Y_{s}\geq 0\,\forall s\in \left[ 0,t\right] ,Y_{t}=b\right)
\end{equation*}%
(see Section 6.1 in \cite{CC2013} for more detail).

In the Brownian case $\alpha =2,\rho =1/2,$ when $\mathcal{Y}$ is a standard
Brownian motion, $\mathbf{P}_{a,b}^{\uparrow ,1}\left( \cdot \right) $ the
law of the so-called normalized Brownian excursion.

For every $N\in \mathbb{N}$ we define a rescaling map $\phi _{N}:\Omega
_{N}^{RW}\rightarrow \Omega _{1}$ by the relation
\begin{equation*}
\left( \phi _{N}(\mathcal{S})\right) (t):=\frac{S_{\left[ Nt\right] }}{a_{N}}%
,
\end{equation*}%
where $\left\{ a_{N}\right\} _{N\in \mathbb{N}}$ is the norming sequence
appearing in Hypothesis A1.

Using this definition we denote by $\mathbb{P}_{x,y}^{\uparrow ,N}\circ \phi
_{N}^{-1}$ the law on $\Omega _{1}:=D\left( [0,1],\mathbb{R}\right) $ given
by the push-forward of $\mathbb{P}_{x,y}^{\uparrow ,N}$ through $\phi _{N}$.

The following important statement was established in \cite{CC2013}.

\begin{theorem}
\label{T_CarCha13} Let $a,b\in \lbrack 0,\infty )$ and let $\left\{
x_{N},N\in \mathbb{N}\right\} $ , $\left\{ y_{N},N\in \mathbb{N}\right\} $
be two non-negative sequences such that $x_{N}/a_{N}\rightarrow a$ and \ $%
y_{N}/a_{N}\rightarrow b$ (in the $(h;c)$-lattice case, assume that $\left(
y_{N}-x_{N}\right) \in Nc+h\mathbb{Z}$ for all $N\in \mathbb{N}$). If
Hypothesis A1 and A2 are satisfied, then, as $N\rightarrow \infty $%
\begin{equation*}
\mathbb{P}_{x,y}^{\uparrow ,N}\circ \phi _{N}^{-1}\Longrightarrow \mathbf{P}%
_{a,b}^{\uparrow ,1}.
\end{equation*}
\end{theorem}

It follows from this theorem that if $a=b=0$ then in the lattice case, for
any function $m=m(N)$ such that $m=o(N)$ as $N\rightarrow \infty $
\begin{equation*}
\lim_{N\rightarrow \infty }\mathbb{P}_{x_{N}}\left( \frac{S_{N-m}}{a_{N}}%
\geq z|L_{N}\geq 0,S_{N}=y_{N}\right) =0
\end{equation*}%
for any $z>0$. Therefore, $S_{N-m}/a_{N}\rightarrow 0$ in probability as $%
N\rightarrow \infty $. Thus, Theorem \ref{T_CarCha13} gives practically no
infromation about the distribution of the random variable $S_{N-m}$ in this
case.

In this section we present the right scaling for $S_{N-m}$ to converge to a
proper non-degenerate random variable depending on the growth rate of $y_{N}$%
. Moreover, to include the absolutely continuous case, we show that%
\begin{equation*}
\lim_{N\rightarrow \infty }\mathbb{P}_{x_{N}}\left( \frac{S_{N-m}-S_{N}}{%
a_{m}}\leq z|L_{N}\geq 0,S_{N}\leq y_{N}\right) =A(z),
\end{equation*}%
where\ $A(z)$ is a proper non-degenerate distribution, whose form is
different for the cases $y_{N}/a_{m}\rightarrow 0,$ $y_{N}/a_{m}\rightarrow
T\in (0,\infty ),$ $y_{N}/a_{m}\rightarrow \infty $. To achieve the goal we
need more notation.

Set $S_{0}:=0,$ $\tau _{k}^{\pm }:=0,$ and denote%
\begin{equation*}
\tau _{k}^{-}:=\inf \left\{ n>\tau _{k-1}^{-}:S_{n}\leq S_{\tau
_{k-1}^{-}}\right\}
\end{equation*}%
for weak descending ladder variables and%
\begin{equation*}
\tau _{k}^{+}:=\inf \left\{ n>\tau _{k-1}^{+}:S_{n}\geq S_{\tau
_{k-1}^{+}}\right\}
\end{equation*}%
for weak ascending ladder variables. \ Put
\begin{equation*}
H_{k}^{\pm }:=\pm S_{\tau _{k}^{\pm }}
\end{equation*}%
and denote%
\begin{eqnarray*}
\zeta &=&\mathbb{P}\left( H_{1}^{+}=0\right) =\sum_{n=1}^{\infty }\mathbb{P}%
\left( S_{1}<0,...,S_{n-1}<0,S_{n}=0\right) \\
&=&\sum_{n=1}^{\infty }\mathbb{P}\left( S_{1}>0,...,S_{n-1}>0,S_{n}=0\right)
=\mathbb{P}\left( H_{1}^{-}=0\right) \in (0,1),
\end{eqnarray*}%
where to justify the third equality it is necessary to use the fact that
\begin{equation*}
\left\{ S_{n}-S_{n-k},k=0,1,...,n\right\} \overset{d}{=}\left\{
S_{k},k=0,1,...,n\right\} .
\end{equation*}

For $x\geq 0$ introduce renewal functions%
\begin{equation*}
V^{\pm }(x)=\sum_{k=0}^{\infty }\mathbb{P}\left( H_{k}^{\pm }\leq x\right)
=\sum_{k=0}^{\infty }\sum_{n=0}^{\infty }\mathbb{P}\left( \tau _{k}^{\pm
}=n,\pm S_{n}\leq x\right) .\quad
\end{equation*}%
Observe that $V^{\pm }(x)$ is a non-decreasing, right-continuous and%
\begin{equation}
V^{\pm }(0)=\sum_{k=0}^{\infty }\mathbb{P}\left( H_{k}^{\pm }=0\right) =%
\frac{1}{1-\zeta }.  \label{V_zero}
\end{equation}%
It is known (see, for instance, \cite{Rog1971}, \cite{Sin57}) that
\begin{equation}
\mathbb{P}\left( \tau _{1}^{+}>n\right) \in R_{-\rho },\quad V^{+}(x)\in
R_{\alpha \rho },  \label{Regular1}
\end{equation}%
and%
\begin{equation}
\mathbb{P}\left( \tau _{1}^{-}>n\right) \in R_{-(1-\rho )},\quad V^{-}(x)\in
R_{\alpha (1-\rho )}.  \label{Regular2}
\end{equation}

Note that, by (15) and (31) in \cite{VW09} and (3.18) in \cite{CC2013} there
are positive constants $\hat{C}$, $\mathcal{C}^{+}$ and \ $\mathcal{C}^{-}$
such that
\begin{equation*}
n\mathbb{P}\left( \tau _{1}^{-}>n\right) \mathbb{P}\left( \tau
_{1}^{+}>n\right) \sim \hat{C},
\end{equation*}%
\begin{equation*}
V^{+}(a_{n})\sim \frac{\mathcal{C}^{+}}{1-\zeta }n\mathbb{P}\left( \tau
_{1}^{-}>n\right) ,\quad V^{-}(a_{n})\sim \frac{\mathcal{C}^{-}}{1-\zeta }n%
\mathbb{P}\left( \tau _{1}^{+}>n\right) ,
\end{equation*}%
and, therefore, as $n\rightarrow \infty $%
\begin{equation}
\mathbb{P}\left( \tau _{1}^{+}>n\right) V^{+}(a_{n})\sim C^{\ast }:=\frac{%
\mathcal{C}^{+}\hat{C}}{1-\zeta }\in (0,\infty ),  \label{Product1}
\end{equation}%
\begin{equation}
\mathbb{P}\left( \tau _{1}^{-}>n\right) V^{-}(a_{n})\sim C^{\ast \ast }:=%
\frac{\mathcal{C}^{-}\hat{C}}{1-\zeta }\in (0,\infty ),  \label{Product11}
\end{equation}%
and%
\begin{equation}
V^{+}(a_{n})V^{-}(a_{n})\sim nC^{\ast \ast \ast }:=\frac{\hat{C}\mathcal{C}%
^{+}\mathcal{C}^{-}}{\left( 1-\zeta \right) ^{2}}n.\quad  \label{Product2}
\end{equation}

We also introduce for $x>0$ left-continuous renewal functions
\begin{equation*}
\underline{{V}}^{\pm }(x):=\sum_{k=0}^{\infty }\mathbb{P}\left( H_{k}^{\pm
}<x\right) =\sum_{k=0}^{\infty }\sum_{n=0}^{\infty }\mathbb{P}\left( \tau
_{k}^{\pm }=n,\pm S_{n}<x\right) .
\end{equation*}%
If the distribution of $X_{1}$ is absolutely continuous, then%
\begin{equation*}
\underline{{V}}^{\pm }(x)=V^{\pm }(x).
\end{equation*}

We will use this fact many times referring to the results established in
\cite{VW09} and \cite{Don12}.

Below we deal with various functions of the form $D(x_{N},y_{N})$ which
depend on real parameters $x_{N}$ and $y_{N}$. Given a fixed positive
sequence $a_{N},n=1,2,...,$ we write

"$D(x_{N},y_{N})=o(1)$ or $D(x_{N},y_{N})\sim c>0$ uniformly in $%
x_{N}=o(a_{N})$, $y_{N}=o(a_{m})$ as $\min (m,N)\rightarrow \infty $"

as a shorthand for

"$D(x_{N},y_{N})=o(1)$ or $D(x_{N},y_{N})\sim c>0$ uniformly in $x_{N}\in
\lbrack 0,\delta _{N}a_{N}]$, $y_{N}\in \lbrack 0,\delta _{N}a_{m}]$ for any
sequence $\delta _{N}\rightarrow 0$ as $\min (m,N)\rightarrow \infty .$"

To simplify the subsequent presentation we recall some conditional local
limit theorems for the lattice case. They are borrowed from Section 4 in
\cite{CC2013}.

Let $g(\cdot )$ be the density of $X_{1}$ and $g^{+}(\cdot )$ be the density
of the time-one marginal distribution of the meander \cite{Chau97} of the
Levy process $\mathcal{Y}$, which may be written informally as
\begin{equation}
g^{+}(x)dx:=\mathbf{P}_{0}\left( Y_{1}\in dx|\inf_{0\leq s\leq 1}Y_{s}\geq
0\right)  \label{Meander1}
\end{equation}%
(see Lemma 4 in \cite{CD2010} for more detail). Analogously, let $%
g^{-}(\cdot )$ be the density of the time-one marginal distribution of the
meander of $-\mathcal{Y}$. Set, finally for $a,b\in \lbrack 0,\infty )$%
\begin{equation*}
C(a,b):=\mathbf{P}_{a}\left( \inf_{0\leq s\leq 1}Y_{s}\geq 0|Y_{1}=b\right) .
\end{equation*}

\begin{lemma}
\label{PropCC}Let Hypothesis A1 and Hypothesis A2 ($(h;c)$-lattice case) be
valid. Then the following relations hold as $n\rightarrow \infty $, for $%
x,y\geq 0$ with $\left( y-x\right) \in nc+h\mathbb{Z}:$

1) uniformly in $x=o(a_{n})$ and $y\geq 0$%
\begin{equation}
g_{n}^{+}\left( x,y\right) =\frac{h\mathbb{P}\left( \tau _{1}^{-}>n\right) }{%
a_{n}}V^{-}(x)\left( g^{+}\left( \frac{y}{a_{n}}\right) +o(1)\right) ;
\label{LocalXsmall}
\end{equation}
2) uniformly in $y=o(a_{n})$ and $x\geq 0$%
\begin{equation}
g_{n}^{+}\left( x,y\right) =\frac{h\mathbb{P}\left( \tau _{1}^{+}>n\right) }{%
a_{n}}V^{+}(y)\left( g^{-}\left( \frac{x}{a_{n}}\right) +o(1)\right) ;
\label{LocalYsmall}
\end{equation}%
3) uniformly in $y=o(a_{n})$ and $x=o(a_{n})$%
\begin{equation}
g_{n}^{+}\left( x,y\right) =h(1-\zeta )\frac{g(0)}{na_{n}}%
V^{-}(x)V^{+}(y)\left( 1+o(1)\right) ;  \label{LocalXYsmall}
\end{equation}%
4) for any $D>1$ uniformly in $x,y\in \left( D^{-1}a_{n},Da_{n}\right) $
\begin{equation}
g_{n}^{+}\left( x,y\right) =\frac{h}{a_{n}}g\left( \frac{y-x}{a_{n}}\right)
C\left( \frac{x}{a_{n}},\frac{y}{a_{n}}\right) \left( 1+o(1)\right) .
\label{LocalXYbig}
\end{equation}
\end{lemma}

We mention that (\ref{LocalXYbig}) is a consequence of Liggett's invariance
principle for bridges \cite{Lig68} and Gnedenko's local limit theorem \cite%
{GK54}\textbf{.}

The renewal functions we have defined above are constructed by weak ladder
variables. Strict ladder variables $\left\{ \hat{\tau}_{k}^{\pm },k\geq
0\right\} $ and $\left\{ \hat{H}_{k}^{\pm },k\geq 0\right\} $ defined as $%
\hat{\tau}_{0}^{\pm }:=0,$\ $\hat{H}_{0}^{\pm }:=0$ and, for $k\geq 1$
\begin{equation*}
\hat{\tau}_{k}^{\pm }:=\inf \left\{ n>\hat{\tau}_{k-1}^{\pm }:\pm S_{n}>\pm
S_{\hat{\tau}_{k-1}^{\pm }}\right\} ,\quad \hat{H}_{k}^{\pm }:=\pm S_{\tau
_{k}^{\pm }},
\end{equation*}
are also often used in studying properties of random walks and branching
processes in random environment (see, for instance, \cite{agkv}, \cite%
{ABGV2011},\cite{Don12}, \cite{VW09},\cite{VD2022}).

The sequences $\left\{ \hat{H}_{k}^{\pm },k\geq 0\right\} $ generate the
renewal functions
\begin{eqnarray*}
\hat{V}^{\pm }(x) &:&=\sum_{k=0}^{\infty }\mathbb{P}\left( \hat{H}_{k}^{\pm
}\leq x\right) =\sum_{k=0}^{\infty }\sum_{n=0}^{\infty }\mathbb{P}\left(
\hat{\tau}_{k}^{\pm }=n,\pm S_{n}\leq x\right) , \\
\underline{\hat{V}}^{\pm }(x) &:&=\sum_{k=0}^{\infty }\mathbb{P}\left( \hat{H%
}_{k}^{\pm }<x\right) =\sum_{k=0}^{\infty }\sum_{n=0}^{\infty }\mathbb{P}%
\left( \hat{\tau}_{k}^{\pm }=n,\pm S_{n}<x\right) .
\end{eqnarray*}%
It is known (see \cite[equation (1.13), Section XII.1]{Fel}) that%
\begin{equation}
\hat{V}^{\pm }(x)=\left( 1-\zeta \right) V^{\pm }(x),\quad \underline{\hat{V}%
}^{\pm }(x)=\left( 1-\zeta \right) \underline{{V}}^{\pm }(x).
\label{RenewRelation}
\end{equation}%
Besides (see \cite[Appendix 3]{CC2013})%
\begin{equation*}
\mathbb{P}\left( \hat{\tau}_{1}^{-}>n\right) \sim \frac{1}{1-\zeta }\mathbb{P%
}\left( \tau _{1}^{-}>n\right)
\end{equation*}%
and, by symmetry,%
\begin{equation}
\mathbb{P}\left( \hat{\tau}_{1}^{+}>n\right) \sim \frac{1}{1-\zeta }\mathbb{P%
}\left( \tau _{1}^{+}>n\right) .  \label{LadderTail2}
\end{equation}

The results analogous to Lemma \ref{PropCC} were obtained in \cite{CC2013}
for the absolutely continuous case as well.

\begin{lemma}
\label{PropCCCont}Let Hypothesis A1 and Hypothesis A2 (absolutely continuous
case) be valid. Then the following relations hold as $n\rightarrow \infty $,
for $x,y\geq 0$

1) uniformly in $x=o(a_{n})$ and $y\geq 0$%
\begin{equation}
f_{n}^{+}\left( x,y\right) =\frac{\mathbb{P}\left( \tau _{1}^{-}>n\right) }{%
a_{n}}V^{-}(x)\left( g^{+}\left( \frac{y}{a_{n}}\right) +o(1)\right) ;
\label{LocXsmallCont}
\end{equation}%
2) uniformly in $y=o(a_{n})$ and $x\geq 0$%
\begin{equation}
f_{n}^{+}\left( x,y\right) =\frac{\mathbb{P}\left( \tau _{1}^{+}>n\right) }{%
a_{n}}V^{+}(y)\left( g^{-}\left( \frac{x}{a_{n}}\right) +o(1)\right) ;
\label{AbsolNew}
\end{equation}%
3) uniformly in $y=o(a_{n})$ and $x=o(a_{n})$%
\begin{equation}
f_{n}^{+}\left( x,y\right) =\frac{g(0)}{na_{n}}V^{-}(x)V^{+}(y)\left(
1+o(1)\right) ;  \label{Don00}
\end{equation}%
4) for any $D>1$ uniformly in $x,y\in \left( D^{-1}a_{n},Da_{n}\right) $
\begin{equation}
f_{n}^{+}\left( x,y\right) =\frac{1}{a_{n}}g\left( \frac{y-x}{a_{n}}\right)
C\left( \frac{x}{a_{n}},\frac{y}{a_{n}}\right) \left( 1+o(1)\right) .
\label{Don12}
\end{equation}
\end{lemma}

In fact, the asymptotic representation (\ref{AbsolNew}) was not mentioned in
\cite{CC2013}. However, it can be deduce from (\ref{LocXsmallCont}) and (\ref%
{ContDistrib}) by symmetry arguments. Namely, considering the random walk $%
-S $ instead of $S$, exchanging $x$ with $y$ and every $"+"$ quantity with
the corresponding $"-"$ one. Relation (\ref{Don12}) is a consequence of
Liggett's invariance principle for bridges \cite{Lig68} and Stone's local
limit theorem \cite{Ston65}\textbf{.}

We now clarify the meaning of constant $C^{\ast }$ in (\ref{Product1}).

\begin{lemma}
\label{L_cStar}Let Hypothesis A1 be valid. Then%
\begin{equation*}
\frac{1}{C^{\ast }}=\int_{0}^{\infty }w^{\alpha \rho }g^{-}(w)dw.
\end{equation*}
\end{lemma}

\begin{remark}
Note that the value of $C^{\ast }$ is not universal in the following sense.
If we change the norming sequence by setting $\bar{a}_{n}=ca_{n}$ with $c>0$
then $S_{n}/\bar{a}_{n}\Longrightarrow \bar{Y}_{1}:=Y_{1}/c$ and, therefore,
\begin{eqnarray*}
\bar{g}^{-}(w)dw &=&\mathbf{P}_{0}\left( -\bar{Y}_{1}\in dw|\inf_{0\leq
s\leq 1}(-\bar{Y}_{s})\geq 0\right) \\
&=&\mathbf{P}_{0}\left( -Y_{1}/c\in dw|\inf_{0\leq s\leq 1}(-Y_{s})\geq
0\right) =cg^{-}(cw)dw.
\end{eqnarray*}%
Thus,%
\begin{eqnarray*}
\int_{0}^{\infty }w^{\alpha \rho }\bar{g}^{-}(w)dw &=&c\int_{0}^{\infty
}w^{\alpha \rho }g^{-}(cw)dw \\
&=&\frac{1}{c^{\alpha \rho }}\int_{0}^{\infty }w^{\alpha \rho }g^{-}(w)dw=%
\frac{1}{c^{\alpha \rho }C^{\ast }},
\end{eqnarray*}%
meaning that $C^{\ast }$ should be replaced by $c^{\alpha \rho }C^{\ast }$.
\end{remark}

\textbf{Proof of Lemma \ref{L_cStar}}. By (\ref{Product11}), (\ref%
{RenewRelation}) and (\ref{LadderTail2})
\begin{equation*}
C^{\ast \ast }=\lim_{n\rightarrow \infty }V^{-}(a_{n})\mathbb{P}\left( \tau
_{1}^{-}>n\right) =\lim_{n\rightarrow \infty }\hat{V}^{-}(a_{n})\mathbb{P}%
\left( \hat{\tau}_{1}^{-}>n\right) .
\end{equation*}%
Let, further $U(w):=C^{\ast \ast }w^{\alpha (1-\rho )},w\geq 0$. It follows
from the proof of Theorem 1.1 in \cite{CC2008} (see the arguments between
formulas (3.11) and (3.12) and definition (3.1)) that%
\begin{equation*}
\mathbf{E}\left[ U(Y_{1})|\inf_{0\leq s\leq 1}Y_{s}\geq 0\right] =C^{\ast
\ast }\int_{0}^{\infty }w^{\alpha (1-\rho )}g^{+}(w)dw=1.
\end{equation*}%
Hence, by considering the time-one marginal distribution of the meander of $-%
\mathcal{Y}$ and observing that the positivity parameter of $-\mathcal{Y}$
is $1-\rho $ we obtain by symmetry arguments for the constant%
\begin{equation*}
C^{\ast }=\lim_{n\rightarrow \infty }V^{+}(a_{n})\mathbb{P}\left( \tau
_{1}^{+}>n\right) =\lim_{n\rightarrow \infty }\hat{V}^{+}(a_{n})\mathbb{P}%
\left( \hat{\tau}_{1}^{+}>n\right)
\end{equation*}%
and the function $\bar{U}(w):=C^{\ast }w^{\alpha \rho }$ that
\begin{equation*}
\mathbf{E}\left[ \bar{U}(-Y_{1})|\inf_{0\leq s\leq 1}(-Y_{s})\geq 0\right]
=C^{\ast }\int_{0}^{\infty }w^{\alpha \rho }g^{-}(w)dw=1.
\end{equation*}

Lemma \ref{L_cStar} is proved

Denote
\begin{equation}
b_{n}:=\frac{1}{na_{n}}=\frac{1}{n^{1+1/\alpha }\ell (n)}.  \label{Def-b}
\end{equation}

Using relations (\ref{RenewRelation}) it is easy to deduce from Proposition
2.3 in \cite{ABGV2011} the following statement which complements Lemmas \ref%
{PropCC} and \ref{PropCCCont}.

\begin{lemma}
\label{Prop_4h}Let Hypothesis A1 hold. There is a number $C>0$ such that,
uniformly for all $x,y\geq 0$ and all $n$%
\begin{equation}
\mathbb{P}_{x}\left( y-1\leq S_{n}<y,L_{n}\geq 0\right) \leq
Cb_{n}V^{-}(x)V^{+}(y)  \label{Rough0}
\end{equation}%
implying%
\begin{equation}
\mathbb{P}_{x}\left( S_{n}<y,L_{n}\geq 0\right) \leq
Cb_{n}V^{-}(x)\sum_{z=0}^{y}V^{+}(z)  \label{Rough1}
\end{equation}%
for the $(h;0)$-lattice case and
\begin{equation}
\mathbb{P}_{x}\left( S_{n}<y,L_{n}\geq 0\right) \leq
Cb_{n}V^{-}(x)\int_{0}^{y}V^{+}(z)dz  \label{Rough11}
\end{equation}%
for the absolutely continuous case.
\end{lemma}

\textbf{Proof}. It was shown in Proposition 2.3 in \cite{ABGV2011} that if
Hypothesis A1 holds, then there is a number $C>0$ such that, uniformly for
all $x,y\geq 0$ and all $n$%
\begin{equation*}
\mathbb{P}_{x}\left( y-1\leq S_{n}<y,L_{n}\geq 0\right) \leq Cb_{n}\hat{V}%
^{-}(x)\underline{\hat{V}}^{+}(y).
\end{equation*}%
Now (\ref{Rough0}) follows from (\ref{RenewRelation}) and the inequalities
\begin{equation*}
\hat{V}^{-}(x)\leq V^{-}(x),\text{ \ \ \ }\underline{\hat{V}}^{+}(y)\leq
V^{+}(y).
\end{equation*}%
The second and third statements of the proposition follow by summation and
integration, respectively.

The proposition is proved.

\section{Limit theorems for lattice case\label{s3}}

Following \cite{CC2013} we consider in the sequel only the $(1;0)$-lattice
case of Hypothesis A2, i.e. we assume that the law of $X_{1}$ is supported
on $\mathbb{Z}$ and is aperiodic. The general $(h;c)$-lattice case requires
a complicated notation and, in other respects, needs not additional
arguments.

Our aim is to study the asymptotic behavior of the probabilities
\begin{equation*}
\mathbb{P}_{x_{N}}\left( S_{N-m}\leq za_{m},L_{N}\geq 0,S_{N}=y_{N}\right)
\end{equation*}
given $\max (x_{N},y_{N})/a_{N}\rightarrow 0$ and $m=o(N)$ as $N\rightarrow
\infty $. We consider three cases: $y_{N}/a_{m}\rightarrow
0,y_{N}/a_{m}\rightarrow T\in (0,\infty ),$ and $y_{N}/a_{m}\rightarrow
\infty $.

\subsection{The case $y_{N}/a_{m}\rightarrow 0\label{s3.1}$}

\begin{lemma}
\label{L_Ysmall} Assume that Hypothesis A1 and A2 are valid, and the
distribution of $X_{1}$ is $(1;0)$-lattice. If $\min (m,N)\rightarrow \infty
$ in such a way that $m=o(N)$, then, for any $z\in (0,\infty ),$
\begin{equation*}
\mathbb{P}_{x_{N}}\left( S_{N-m}\leq za_{m}|L_{N}\geq 0,S_{N}=y_{N}\right)
=A_{1}(z)(1+o(1))
\end{equation*}%
uniformly in $x_{N}=o(a_{N})$ and $y_{N}=o(a_{m})$, where (recall (\ref%
{Product1}) and Lemma \ref{L_cStar})%
\begin{equation}
A_{1}(z):=C^{\ast }\int_{0}^{z}w^{\alpha \rho }g^{-}\left( w\right) dw,\quad
z\in \lbrack 0,\infty ),  \label{DefA1(z)}
\end{equation}%
is a proper distribution.
\end{lemma}

\textbf{Proof. }According to (\ref{LocalXYsmall})
\begin{equation}
\mathbb{P}_{x_{N}}\left( L_{N}\geq 0,S_{N}=j\right) \sim (1-\zeta
)b_{N}g(0)V^{-}(x_{N})V^{+}(j)  \label{BasicAsympt}
\end{equation}%
uniformly in $x_{N}=o\left( a_{N}\right) $ and $j=o\left( a_{N}\right) $
implying%
\begin{equation}
\mathbb{P}_{x_{N}}\left( L_{N}\geq 0,S_{N}\leq y_{N}\right) \sim (1-\zeta
)b_{N}g(0)V^{-}(x_{N})\sum_{j=0}^{y_{N}}V^{+}(j)  \label{BasicAsympt2}
\end{equation}%
uniformly in $x_{N}=o\left( a_{N}\right) $ and $y_{N}=o\left( a_{N}\right) $%
. Hence, in view of (\ref{V_zero})%
\begin{equation}
\mathbb{P}_{x_{N}}\left( L_{N}\geq 0,S_{N}\leq y_{N}\right) \sim (1-\zeta
)V^{-}(x_{N})\mathbb{P}\left( L_{N}\geq 0,S_{N}\leq y_{N}\right)
\label{BasicAsympt3}
\end{equation}%
uniformly in $x_{N}=o\left( a_{N}\right) $ and $y_{N}=o\left( a_{N}\right) $.

We now fix $z>\varepsilon \in (0,1)$ and investigate the asymptotic behavior
of the probability
\begin{equation}
\mathbb{P}_{x_{N}}\left( S_{N-m}\leq za_{m},L_{N}\geq 0,S_{n}=y_{N}\right)
=R_{1}(\varepsilon ,m,N)+R_{2}(\varepsilon ,m,N),  \label{Decompos}
\end{equation}%
where
\begin{eqnarray*}
R_{1}(\varepsilon ,m,N)&:=&\sum_{0\leq k\leq \varepsilon a_{m}}\mathbb{P}%
_{x_{N}}\left( S_{N-m}=k,L_{N-m}\geq 0\right) \mathbb{P}_{k}\left(
S_{m}=y_{N},L_{m}\geq 0\right) , \\
R_{2}(\varepsilon ,m,N)&:=&\sum_{\varepsilon a_{m}<k\leq za_{m}}\mathbb{P}%
_{x_{N}}\left( S_{N-m}=k,L_{N-m}\geq 0\right) \mathbb{P}_{k}\left(
S_{m}=y_{N},L_{m}\geq 0\right) .
\end{eqnarray*}

Using Lemma \ref{Prop_4h} we conclude that
\begin{eqnarray*}
R_{1}(\varepsilon ,m,N) &\leq &Cb_{m}V^{+}(y_{N})\sum_{0\leq k\leq
\varepsilon a_{m}}\mathbb{P}_{x_{N}}\left( S_{N-m}=k,L_{N-m}\geq 0\right)
V^{-}(k) \\
&\leq &Cb_{m}V^{+}(y_{N})V^{-}(\varepsilon a_{m})\mathbb{P}_{x_{N}}\left(
S_{N-m}\leq \varepsilon a_{m},L_{N-m}\geq 0\right) \\
&\leq &C_{1}b_{m}b_{N-m}V^{+}(y_{N})V^{-}(x_{N})V^{-}(\varepsilon
a_{m})\sum_{0\leq k\leq \varepsilon a_{m}}V^{+}(k) \\
&\leq &\varepsilon
C_{1}a_{m}b_{m}b_{N-m}V^{+}(y_{N})V^{-}(x_{N})V^{-}(\varepsilon
a_{m})V^{+}(\varepsilon a_{m}).
\end{eqnarray*}%
Applying (\ref{Product2}) and (\ref{Def-b}) and using the equivalence $%
b_{N-m}\sim b_{N}$ for sufficiently large $N$ and $m=o(N)$ we get
\begin{eqnarray*}
&&a_{m}b_{m}b_{N-m}V^{+}(y_{N})V^{-}(x_{N})V^{-}(\varepsilon
a_{m})V^{+}(\varepsilon a_{m}) \\
&&\qquad\leq Ca_{m}b_{m}b_{N}V^{+}(y_{N})V^{-}(x_{N})V^{-}(a_{m})V^{+}(a_{m})
\\
&&\qquad\leq 2CC^{\ast \ast \ast }a_{m}b_{m}b_{N}mV^{+}(y_{N})V^{-}(x_{N}) \\
&&\qquad\quad=2CC^{\ast \ast \ast }b_{N}V^{+}(y_{N})V^{-}(x_{N}).
\end{eqnarray*}%
Thus,%
\begin{equation}
R_{1}(\varepsilon ,m,N)\leq \varepsilon C_{1}b_{N}V^{+}(y_{N})V^{-}(x_{N}).
\label{R_1}
\end{equation}

To evaluate $R_{2}(\varepsilon,m,N)$ we first observe that according to
Lemma \ref{PropCC}
\begin{equation}
\mathbb{P}_{x_{N}}\left( S_{N-m}=k,L_{N-m}\geq 0\right) \sim \frac{(1-\zeta
)g(0)}{\left( N-m\right) a_{N-m}}V^{-}(x_{N})V^{+}(k)  \label{R_2}
\end{equation}%
uniformly in $\varepsilon a_{m}\leq k\leq za_{m}$ and $%
x_{N}=o(a_{N-m})=o(a_{N}),$ and%
\begin{equation}
\mathbb{P}_{k}\left( S_{m}=y_{N},L_{m}\geq 0\right) =\frac{\mathbb{P}\left(
\tau _{1}^{+}>m\right) }{a_{m}}V^{+}(y_{N})\left( g^{-}\left( \frac{k}{a_{m}}%
\right) +o(1)\right)  \label{R_3}
\end{equation}%
uniformly in $\varepsilon a_{m}\leq k\leq za_{m}$ and $y_{N}=o(a_{m})$.
Thus, uniformly in $\varepsilon a_{m}\leq k\leq za_{m}$ and $y_{N}=o(a_{m})$
and $x_{N}=o(a_{N})$
\begin{eqnarray*}
R_{2}(\varepsilon ,m,N) &\sim &\frac{(1-\zeta )g(0)}{Na_{N}}\frac{\mathbb{P}%
\left( \tau _{1}^{+}>m\right) }{a_{m}}V^{-}(x_{N})V^{+}(y_{N}) \\
&&\times \sum_{\varepsilon a_{m}<k\leq za_{m}}V^{+}(k)g^{-}\left( \frac{k}{%
a_{m}}\right) .
\end{eqnarray*}%
Using (\ref{Product1}) we conclude that%
\begin{eqnarray*}
R_{2}(\varepsilon ,m,N) &\sim &\frac{C^{\ast }(1-\zeta )g(0)}{Na_{N}}%
V^{-}(x_{N})V^{+}(y_{N}) \\
&&\times \sum_{\varepsilon a_{m}<k\leq za_{m}}\frac{V^{+}(k)}{V^{+}(a_{m})}%
g^{-}\left( \frac{k}{a_{m}}\right) \frac{1}{a_{m}}.
\end{eqnarray*}%
It follows from (\ref{Regular1}) and properties of regularly varying
functions (see \cite{sen76}) that, as $m\rightarrow \infty $
\begin{equation*}
\frac{V^{+}(wa_{m})}{V^{+}(a_{m})}\rightarrow w^{\alpha \rho }
\end{equation*}%
uniformly in $\varepsilon \leq w\leq z$. Hence we get, as $m\rightarrow
\infty $%
\begin{equation*}
\sum_{\varepsilon a_{m}<k\leq za_{m}}\frac{V^{+}(k)}{V^{+}(a_{m})}%
g^{-}\left( \frac{k}{a_{m}}\right) \frac{1}{a_{m}}\sim \int_{\varepsilon
}^{z}w^{\alpha \rho }g^{-}\left( w\right) dw.
\end{equation*}%
As a result,%
\begin{equation*}
R_{2}(\varepsilon ,m,N)\sim C^{\ast }(1-\zeta
)g(0)b_{N}V^{-}(x_{N})V^{+}(y_{N})\int_{\varepsilon }^{z}w^{\alpha \rho
}g^{-}\left( w\right) dw
\end{equation*}%
as $N\rightarrow \infty $ uniformly in $y_{N}=o(a_{m})$ and $x_{N}=o(a_{N})$%
. Since $\varepsilon >0$ may be selected arbitrary small, we conclude by (%
\ref{R_1}) and (\ref{BasicAsympt}) that
\begin{equation}
\frac{\mathbb{P}_{x_{N}}\left( S_{N-m}\leq za_{m},L_{N}\geq
0,S_{N}=y_{N}\right) }{\mathbb{P}_{x_{N}}\left( L_{N}\geq
0,S_{N}=y_{N}\right) }\sim C^{\ast }\int_{0}^{z}w^{\alpha \rho }g^{-}\left(
w\right) dw=A_{1}(z)  \label{C00}
\end{equation}%
uniformly in $x_{N}=o(a_{N})$ and $y_{N}=o(a_{m})$.

The lemma is proved.

\begin{corollary}
\label{C_Ysmall} Assume that Hypothesis A1 and A2 are valid, and the
distribution of $X_{1}$ is $(1;0)$-lattice. If $\min (m,N)\rightarrow \infty
$ in such a way that $m=o(N)$, then, for any $z\in (0,\infty )$
\begin{equation}
\mathbb{P}_{x_{N}}\left( S_{N-m}\leq za_{m},L_{N}\geq 0,S_{N}\leq
y_{N}\right) \sim A_{1}(z)\mathbb{P}_{x_{N}}\left( L_{N}\geq 0,S_{N}\leq
y_{N}\right)  \label{C11}
\end{equation}%
uniformly in $x_{N}=o(a_{N})$ and $y_{N}=o(a_{m})$.
\end{corollary}

\textbf{Proof}. The desired statement follows from (\ref{C00}), (\ref%
{BasicAsympt}), and (\ref{BasicAsympt2}) by summation.

\subsection{The case $y_{N}\asymp a_{m}\label{s3.2}$}

Let $0<t_{0}<t_{1}<\infty $ be fixed.

\begin{lemma}
\label{L_Yequal}Let Hypothesis A1 and A2 be valid, and the distribution of $%
X_{1}$ be $(1;0)$-lattice. If $\min (m,N)\rightarrow \infty $ in such a way
that $m=o(N)$, then, for any $z\in (0,\infty )$
\begin{equation*}
\mathbb{P}_{x_{N}}\left( S_{N-m}\leq za_{m}|L_{N}\geq 0,S_{N}=j\right) \sim
A_{2}\left( z,\frac{j}{a_{m}}\right)
\end{equation*}%
uniformly in $x_{N}=o(a_{N})$ and $j\in \lbrack t_{0}a_{m},t_{1}a_{m}]$,
where%
\begin{equation*}
A_{2}(z,t):=t^{-\alpha \rho }\int_{0}^{z}w^{\alpha \rho }g\left( t-w\right)
C\left( w,t\right) dw.
\end{equation*}
\end{lemma}

\textbf{Proof}. We again use the decomposition (\ref{Decompos}). The
estimates (\ref{R_1})\ and (\ref{R_2}) remain the same while instead of (\ref%
{R_3}) we have by (\ref{LocalXYbig}) as $m\rightarrow \infty $
\begin{equation}
\mathbb{P}_{k}\left( S_{m}=j,L_{m}\geq 0\right) \sim \frac{1}{a_{m}}g\left(
\frac{j-k}{a_{m}}\right) C\left( \frac{k}{a_{m}},\frac{j}{a_{m}}\right)
\label{R_4}
\end{equation}%
uniformly in $\varepsilon a_{m}\leq k\leq za_{m}$ and $j\in \lbrack
t_{0}a_{m},t_{1}a_{m}]$. Hence we deduce by (\ref{R_2}) and (\ref{R_4}) that
\begin{eqnarray*}
R_{2}(\varepsilon ,m,N) &\sim &\frac{(1-\zeta )g(0)}{\left( N-m\right)
a_{N-m}}\text{$V$}^{-}(x_{N}) \\
&&\times \sum_{\varepsilon a_{m}<k\leq za_{m}}\text{$V$}^{+}(k)\frac{1}{a_{m}%
}g\left( \frac{j-k}{a_{m}}\right) C\left( \frac{k}{a_{m}},\frac{j}{a_{m}}%
\right) \\
&\sim &\frac{(1-\zeta )g(0)}{Na_{N}}\text{$V$}^{-}(x_{N})\text{$V$}%
^{+}(a_{m}) \\
&&\times \sum_{\varepsilon a_{m}<k\leq za_{m}}\frac{\text{$V$}^{+}(k)}{\text{%
$V$}^{+}(a_{m})}\frac{1}{a_{m}}g\left( \frac{j-k}{a_{m}}\right) C\left(
\frac{k}{a_{m}},\frac{j}{a_{m}}\right)
\end{eqnarray*}%
uniformly in $x_{N}=o(a_{N})$ and $j\in \lbrack t_{0}a_{m},t_{1}a_{m}]$.
Now, the same as in Lemma \ref{L_Ysmall} we have
\begin{equation*}
\sum_{\varepsilon a_{m}<k\leq za_{m}}\frac{\text{$V$}^{+}(k)}{\text{$V$}%
^{+}(a_{m})}\frac{1}{a_{m}}g\left( \frac{j-k}{a_{m}}\right) C\left( \frac{k}{%
a_{m}},\frac{j}{a_{m}}\right) \sim \int_{\varepsilon }^{z}w^{\alpha \rho
}g\left( \frac{j}{a_{m}}-w\right) C\left( w,\frac{j}{a_{m}}\right) dw
\end{equation*}%
as $m\rightarrow \infty $. Thus,%
\begin{eqnarray}
R_{2}(\varepsilon ,m,N) &\sim &(1-\zeta )g(0)b_{N}\text{$V$}^{-}(x_{N})\text{%
$V$}^{+}(a_{m})  \notag \\
&&\times \int_{\varepsilon }^{z}w^{\alpha \rho }g\left( \frac{j}{a_{m}}%
-w\right) C\left( w,\frac{j}{a_{m}}\right) dw.  \label{R_34}
\end{eqnarray}%
Since $\varepsilon >0$ may be selected arbitrary small, and
\begin{equation*}
V^{+}(a_{m})=V^{+}\left( \frac{a_{m}}{j}j\right) \sim \left( \frac{a_{m}}{j}%
\right) ^{\alpha \rho }V^{+}\left( j\right)
\end{equation*}%
uniformly in $j\in \lbrack t_{0}a_{m},t_{1}a_{m}],$ we conclude by (\ref{R_1}%
), (\ref{R_34}) and (\ref{BasicAsympt}) that%
\begin{equation*}
\frac{\mathbb{P}_{x_{N}}\left( S_{N-m}\leq za_{m},L_{N}\geq 0,S_{N}=j\right)
}{\mathbb{P}_{x_{N}}\left( L_{N}\geq 0,S_{N}=j\right) }=A_{2}\left( z,\frac{j%
}{a_{m}}\right) (1+o(1))
\end{equation*}%
uniformly in $x_{N}=o(a_{N})$ and $j\in \lbrack t_{0}a_{m},t_{1}a_{m}]$.

Lemma \ref{L_Yequal} is proved.

\begin{corollary}
\label{C_Yequal}If $\min (m,N)\rightarrow \infty ,m=o(N)$ and $y_{N}\sim
Ta_{m},T\in (0,\infty ),$ then, under the conditions of Lemma \ref{L_Yequal}%
, for any $z\in (0,\infty )$%
\begin{equation*}
\mathbb{P}_{x_{N}}\left( S_{N-m}\leq za_{m},L_{N}\geq 0,S_{N}\leq
y_{N}\right) \sim \mathbb{P}_{x_{N}}\left( L_{N}\geq 0,S_{N}\leq
y_{N}\right) B(z,T),
\end{equation*}%
uniformly in $x_{N}=o(a_{N})$, where%
\begin{equation}
B(z,T)=\frac{\alpha \rho +1}{T^{\alpha \rho +1}}\int_{0}^{z}w^{\alpha \rho
}dw\int_{0}^{T}g\left( t-w\right) C\left( w,t\right) dt.  \label{DefB(z,T)}
\end{equation}
\end{corollary}

\textbf{Proof}. In view of (\ref{Rough1}) for any $\varepsilon \in (0,1)$%
\begin{eqnarray}
&&\mathbb{P}_{x_{N}}\left( S_{N-m}\leq za_{m},L_{N}\geq 0,S_{N}\leq
\varepsilon y_{N}\right)  \notag \\
&&\quad \leq \mathbb{P}_{x_{N}}\left( L_{N}\geq 0,S_{N}\leq \varepsilon
y_{N}\right)  \notag \\
&&\qquad \leq C_{1}b_{N}V^{-}(x_{N})\sum_{0\leq j\leq \varepsilon
y_{N}}V^{+}(j).  \label{SmallYEqual1}
\end{eqnarray}%
Further, by Lemma \ref{L_Yequal}, and relations (\ref{LocalXYsmall}) and (%
\ref{Regular1})
\begin{eqnarray}
&&\mathbb{P}_{x_{N}}\left( S_{N-m}\leq za_{m},L_{N}\geq 0,\varepsilon
y_{N}\leq S_{N}\leq y_{N}\right)  \notag \\
&&\quad \sim \sum_{\varepsilon y_{N}<j\leq y_{N}}\mathbb{P}_{x_{N}}\left(
L_{N}\geq 0,S_{N}=j\right) A_{2}(z,j/a_{m})  \notag \\
&&\quad \sim (1-\zeta )g(0)b_{N}V^{-}(x_{N})\sum_{\varepsilon y_{N}<j\leq
y_{N}}V^{+}(j)A_{2}(z,j/a_{m})  \notag \\
&&\quad =(1-\zeta )g(0)b_{N}V^{-}(x_{N})V^{+}(a_{m})\sum_{\varepsilon
Ta_{m}<j\leq Ta_{m}}\frac{V^{+}(j)}{V^{+}(a_{m})}A_{2}(z,j/a_{m})  \notag \\
&&\quad \sim (1-\zeta
)g(0)b_{N}V^{-}(x_{N})V^{+}(a_{m})a_{m}\int_{\varepsilon T}^{T}t^{\alpha
\rho }A_{2}(z,t)dt.  \label{EstimMain}
\end{eqnarray}%
Note that%
\begin{eqnarray*}
\int_{\varepsilon T}^{T}t^{\alpha \rho }A_{2}(z,t)dt &=&\int_{\varepsilon
T}^{T}dt\int_{0}^{z}w^{\alpha \rho }g\left( t-w\right) C\left( w,t\right) dw
\\
&=&\int_{0}^{z}w^{\alpha \rho }dw\int_{\varepsilon T}^{T}g\left( t-w\right)
C\left( w,t\right) dt.
\end{eqnarray*}%
Since
\begin{equation}
V^{+}(a_{m})a_{m}\sim V^{+}(T^{-1}y_{N})T^{-1}y_{N}\sim T^{-\alpha \rho
-1}V^{+}(y_{N})y_{N}  \label{EquivT}
\end{equation}%
by (\ref{Regular1}), we get in view of (\ref{EstimMain})
\begin{eqnarray*}
&&\mathbb{P}_{x_{N}}\left( S_{N-m}\leq za_{m},L_{N}\geq 0,\varepsilon
y_{N}\leq S_{N}\leq y_{N}\right) \\
&\sim &(1-\zeta )g(0)b_{N}V^{-}(x_{N})V^{+}(y_{N})y_{N} \\
&&\times T^{-\alpha \rho -1}\int_{0}^{z}w^{\alpha \rho }dw\int_{\varepsilon
T}^{T}g\left( t-w\right) C\left( w,t\right) dt.
\end{eqnarray*}%
Recalling (\ref{Regular1}) once again we see that%
\begin{equation*}
\left( \alpha \rho +1\right) \sum_{j=0}^{y_{N}}V^{+}(j)\sim y_{N}V^{+}(y_{N})
\end{equation*}%
as $y_{N}\rightarrow \infty $. Hence, using (\ref{SmallYEqual1}), (\ref%
{BasicAsympt2}) and letting $\varepsilon $ to zero we conclude that%
\begin{equation*}
\mathbb{P}_{x_{N}}\left( S_{N-m}\leq za_{m},L_{N}\geq 0,S_{N}\leq
y_{N}\right) \sim \mathbb{P}_{x_{N}}\left( L_{N}\geq 0,S_{N}\leq
y_{N}\right) B(z,T).
\end{equation*}%
The corollary is proved.

\subsection{The case $a_{m}=o(y_{N})$}

In contrast to sections \ref{s3.1} and \ref{s3.2}, where we have
investigated the distribution of the random variable $S_{N-m}$, we analyze
here the distribution of the difference $S_{N-m}-S_{N}$.

\begin{lemma}
\label{L_Ybig}Let Hypothesis A1 and A2 be valid, and the distribution of $%
X_{1}$ be $(1;0)$-lattice. If $\min (m,N)\rightarrow \infty $ and $%
a_{m}=o(y_{N})$ then, for any $z\in (-\infty ,\infty )$
\begin{equation*}
\mathbb{P}_{x_{N}}\left( S_{N-m}-S_{N}\leq za_{m}|L_{N}\geq
0,S_{N}=y_{N}\right) \sim \mathbf{P}\left( Y_{1}\leq z\right)
\end{equation*}%
uniformly in $x_{N}=o(a_{N})$ and $y_{N}=o(a_{N})$, $a_{m}=o(y_{N})$, where $%
Y_{1}$ is defined by~(\ref{LimitLaw}).
\end{lemma}

\textbf{Proof}. We fix a\ large value $M$ and use the decomposition%
\begin{equation*}
\mathbb{P}_{x_{N}}\left( S_{N-m}-y_{N}\leq za_{m},L_{N}\geq
0,S_{N}=y_{N}\right) =R_{3}(M,m,N)+R_{4}(M,m,N),
\end{equation*}%
where%
\begin{eqnarray*}
R_{3}(M,m,N) &:&=\sum_{0\leq k\leq y_{N}-Ma_{m}}\mathbb{P}_{x_{N}}\left(
S_{N-m}=k,L_{N-m}\geq 0\right) \mathbb{P}_{k}\left( S_{m}=y_{N},L_{m}\geq
0\right) , \\
R_{4}(M,m,N) &:&=\sum_{y_{N}-Ma_{m}<k\leq y_{N}+za_{m}}\mathbb{P}%
_{x_{N}}\left( S_{N-m}=k,L_{N-m}\geq 0\right) \mathbb{P}_{k}\left(
S_{m}=y_{N},L_{m}\geq 0\right) .
\end{eqnarray*}%
Observe that by (\ref{LocalXYsmall})
\begin{eqnarray*}
R_{3}(M,m,N) &\leq &Cb_{N-m}\text{$V$}^{-}(x_{N})\sum_{0\leq k\leq
y_{N}-Ma_{m}}\text{$V$}^{+}(k)\mathbb{P}_{k}\left( S_{m}=y_{N},L_{m}\geq
0\right) \\
&\leq &C_{1}b_{N}\text{$V$}^{-}(x_{N})\text{$V$}^{+}(y_{N})\sum_{0\leq k\leq
y_{N}-Ma_{m}}\mathbb{P}_{k}\left( S_{m}=y_{N},L_{m}\geq 0\right) .
\end{eqnarray*}%
Now

\begin{align*}
\sum_{0\leq k\leq y_{N}-Ma_{m}}\mathbb{P}_{k}\left( S_{m}=y_{N},L_{m}\geq
0\right) =\sum_{0\leq k\leq y_{N}-Ma_{m}}\mathbb{P}\left(
S_{m}=y_{N}-k,L_{m}\geq -k\right) \\
\leq\sum_{0\leq k\leq y_{N}-Ma_{m}}\mathbb{P}\left( S_{m}=y_{N}-k\right)
\leq \mathbb{P}\left( S_{m}\geq Ma_{m}\right) .
\end{align*}%
Since $\lim_{M\rightarrow \infty }\lim_{m\rightarrow \infty }\mathbb{P}%
\left( S_{m}\geq Ma_{m}\right) =0$, it follows that%
\begin{equation}
R_{3}(M,m,N)\leq \varepsilon (M)b_{N}\text{$V$}^{-}(x_{N})\text{$V$}%
^{+}(y_{N}),  \label{R_33}
\end{equation}%
where $\varepsilon (M)\downarrow 0$ as $M\rightarrow \infty .$

Further, in view of (\ref{LocalXYsmall}) and (\ref{Regular1})\
\begin{eqnarray*}
R_{4}(M,m,N) &\sim &(1-\zeta )g(0)b_{N-m}\text{$V$}^{-}(x_{N})%
\sum_{y_{N}-Ma_{m}<k\leq y_{N}+za_{m}}\text{$V$}^{+}(k)\mathbb{P}_{k}\left(
S_{m}=y_{N},L_{m}\geq 0\right) \\
&\sim &(1-\zeta )g(0)b_{N}\text{$V$}^{-}(x_{N})\text{$V$}^{+}(y_{N})%
\sum_{y_{N}-Ma_{m}<k\leq y_{N}+za_{m}}\mathbb{P}_{k}\left(
S_{m}=y_{N},L_{m}\geq 0\right)
\end{eqnarray*}%
uniformly in $x_{N}=o(a_{N})$ and $y_{N}=o(a_{N})$. Thus, it remains to
evaluate the sum
\begin{eqnarray*}
&&\sum_{y_{N}-Ma_{m}<k\leq y_{N}+za_{m}}\mathbb{P}_{k}\left(
S_{m}=y_{N},L_{m}\geq 0\right) \\
&&\qquad \qquad =\sum_{y_{N}-Ma_{m}<k\leq y_{N}+za_{m}}\mathbb{P}\left(
S_{m}=y_{N}-k,L_{m}\geq -k\right) \\
&&\qquad \qquad =\sum_{-Ma_{m}<j\leq za_{m}}\mathbb{P}\left(
S_{m}=j,L_{m}\geq j-y_{N}\right) .
\end{eqnarray*}%
Clearly,%
\begin{align*}
\mathbb{P}\left( S_{m}\in \lbrack -Ma_{m},za_{m}],L_{m}\geq
za_{m}-y_{N}\right) & \leq \sum_{-Ma_{m}<j\leq za_{m}}\mathbb{P}\left(
S_{m}=j,L_{m}\geq j-y_{N}\right) \\
& \leq \mathbb{P}\left( S_{m}\in \lbrack -Ma_{m},za_{m}]\right) .
\end{align*}%
Further,%
\begin{eqnarray*}
\mathbb{P}\left( S_{m}\in \lbrack -Ma_{m},za_{m}],L_{m}\geq
za_{m}-y_{N}\right) &\geq &\mathbb{P}\left( S_{m}\in \lbrack
-Ma_{m},za_{m}]\right) \\
&&-\mathbb{P}\left( L_{m}<za_{m}-y_{N}\right) .
\end{eqnarray*}%
Since $y_{N}/a_{m}\rightarrow \infty $, it follows by the invariance
principle for random walks whose increments have distributions belonging
without centering to the domain of attraction of a stable law that
\begin{equation*}
\lim_{N\rightarrow \infty }\mathbb{P}\left( \frac{L_{m}}{a_{m}}<z-\frac{y_{N}%
}{a_{m}}\right) =0.
\end{equation*}%
Observing that $b_{N-m}\sim b_{N}$ if $m=o(N)$ as $N\rightarrow \infty $, we
get
\begin{equation*}
R_{4}(M,m,N)\sim (1-\zeta )g(0)b_{N}\text{$V$}^{-}(x_{N})\text{$V$}%
^{+}(y_{N})\mathbf{P}\left( Y_{1}\in \left[ -M,z\right] \right)
\end{equation*}%
as $N\rightarrow \infty $. Since $M$ may be selected arbitrary large and $%
\varepsilon (M)>0$ in (\ref{R_33})\ arbitrary small, we conclude that%
\begin{align*}
& \mathbb{P}_{x_{N}}\left( S_{N-m}-y_{N}\leq za_{m};L_{N}\geq
0,S_{n}=y_{N}\right) \\
\qquad & \sim (1-\zeta )g(0)b_{N}\text{$V$}^{-}(x_{N})\text{$V$}^{+}(y_{N})%
\mathbf{P}\left( Y_{1}\leq z\right)
\end{align*}%
uniformly in $x_{N}=o(a_{N})$ and $y_{N}=o(a_{N})$. This, in view of (\ref%
{LimitLaw}) implies the statement of the lemma.

\begin{corollary}
\label{C_Ybig}Let Hypothesis A1 and A2 be valid, and the distribution of $%
X_{1}$ be $(1;0)$-lattice. If $m\rightarrow \infty $ and $a_{m}=o(y_{N})$
then, for any $z\in (-\infty ,\infty )$
\begin{equation*}
\mathbb{P}_{x_{N}}\left( S_{N-m}-S_{N}\leq za_{m},L_{N}\geq 0,S_{N}\leq
y_{N}\right) \sim \mathbb{P}_{x_{N}}\left( L_{N}\geq 0,S_{N}\leq
y_{N}\right) \mathbf{P}\left( Y_{1}\leq z\right)
\end{equation*}%
uniformly in $x_{N}=o(a_{N})$ and $y_{N}=o(a_{N}),a_{m}=o(y_{N})$.
\end{corollary}

\textbf{Proof}. In view of (\ref{Rough1}) for any $\varepsilon \in (0,1)$%
\begin{equation*}
\mathbb{P}_{x_{N}}\left( S_{N-m}-S_{N}\leq za_{m};L_{N}\geq 0,S_{N}\leq
\varepsilon y_{N}\right) \leq C_{1}b_{N}V^{-}(x_{N})\sum_{0\leq j\leq
\varepsilon y_{N}}V^{+}(j).
\end{equation*}%
Further, by Lemma \ref{L_Ybig}
\begin{eqnarray*}
&&\mathbb{P}_{x_{N}}\left( S_{N-m}-S_{N}\leq za_{m};L_{N}\geq 0,\varepsilon
y_{N}\leq S_{N}\leq y_{N}\right) \\
&&\qquad \sim \sum_{\varepsilon y_{N}<j\leq y_{N}}\mathbb{P}_{x_{N}}\left(
L_{N}\geq 0,S_{N}=j\right) \mathbf{P}\left( Y_{1}\leq z\right) \\
&&\qquad \sim \mathbb{P}_{x_{N}}\left( L_{N}\geq 0,\varepsilon y_{N}\leq
S_{N}\leq y_{N}\right) \mathbf{P}\left( Y_{1}\leq z\right) .
\end{eqnarray*}%
In view of (\ref{LocalXYsmall})\
\begin{equation*}
\mathbb{P}_{x_{N}}\left( L_{N}\geq 0,S_{N}\leq y_{N}\right) \sim (1-\zeta
)g(0)b_{N}\text{$V$}^{-}(x_{N})\sum_{j=0}^{y_{N}}\text{$V$}^{+}(j).
\end{equation*}%
\ Hence the desired estimate follows by using (\ref{Regular1}) and letting $%
\varepsilon \downarrow 0$.

\section{Limit theorems for absolutely continuous case\label{s4}}

In this section we deal with analogues of Corollaries \ref{C_Ysmall}, \ref%
{C_Yequal} and \ref{C_Ybig} for the absolutely continuous case.

\begin{lemma}
\label{L_TotalCont}Let Hypothesis A1 and Hypothesis A2 (absolutely
continuous case) hold and $\min (m,N)\rightarrow \infty ,m=o(N)$. Then

1) uniformly in $x_{N}=o(a_{N})$ and $y_{N}=o(a_{m})$
\begin{equation*}
\mathbb{P}_{x_{N}}\left( S_{N-m}\leq za_{m}|L_{N}\geq 0,S_{N}\leq
y_{N}\right) \sim A_{1}(z),\ z\in (0,\infty );
\end{equation*}

2) If $y_{N}\sim Ta_{m},T\in (0,\infty )$ then, uniformly in $x_{N}=o(a_{N})$
\begin{equation}
\mathbb{P}_{x_{N}}\left( S_{N-m}\leq za_{m}|L_{N}\geq 0,S_{N}\leq
y_{N}\right) \sim B(z,T),\ z\in (0,\infty );  \label{Con333}
\end{equation}

3) if $a_{m}=o(y_{N})$ then, uniformly in $x_{N}=o(a_{N})$ and $%
y_{N}=o(a_{N})$%
\begin{equation*}
\mathbb{P}_{x_{N}}\left( S_{N-m}-S_{N}\leq za_{m}|L_{N}\geq 0,S_{N}\leq
y_{N}\right) \sim \mathbf{P}(Y_{1}\leq z),\ z\in (-\infty ,\infty ).
\end{equation*}
\end{lemma}

\textbf{Proof}. The proofs of the statements of the lemma have practically
no changes with the proofs of the respective statements of Corollaries \ref%
{C_Ysmall}, \ref{C_Yequal} and \ref{C_Ybig}. For this reason we check only (%
\ref{Con333}).

In view of (\ref{Rough0}) for any $\varepsilon \in (0,1)$%
\begin{eqnarray*}
\mathbb{P}_{x_{N}}\left( S_{N-m}\leq za_{m},L_{N}\geq 0,S_{N}<\varepsilon
y_{N}\right) &\leq &\mathbb{P}_{x_{N}}\left( L_{N}\geq 0,S_{N}<\varepsilon
y_{N}\right) \\
&\leq &C_{1}b_{N-m}V^{-}(x_{N})\int_{0}^{\varepsilon y_{N}}V^{+}(w)dw.
\end{eqnarray*}%
Similarly,%
\begin{eqnarray*}
\mathbb{P}_{x_{N}}\left( S_{N-m}<\varepsilon y_{N},L_{N}\geq 0,S_{N}\leq
y_{N}\right) &\leq &\mathbb{P}_{x_{N}}\left( S_{N-m}<\varepsilon
y_{N},L_{N-m}\geq 0\right) \\
&\leq &C_{1}b_{N-m}V^{-}(x_{N})\int_{0}^{\varepsilon y_{N}}V^{+}(w)dw.
\end{eqnarray*}%
Further,
\begin{eqnarray*}
&&\mathbb{P}_{x_{N}}\left( \varepsilon y_{N}\leq S_{N-m}\leq
za_{m},L_{N}\geq 0,\varepsilon y_{N}\leq S_{N}\leq y_{N}\right) \\
&&\quad=\int_{\varepsilon y_{N}}^{za_{m}}\mathbb{P}_{x_{N}}\left( S_{N-m}\in
dw,L_{N-m}\geq 0\right) \mathbb{P}_{w}\left( L_{m}\geq 0,\varepsilon
y_{N}\leq S_{m}\leq y_{N}\right) .
\end{eqnarray*}%
By (\ref{Don12}) for any $0<t_{0}<t_{1}<\infty $
\begin{eqnarray*}
\mathbb{P}_{w}\left( L_{m}\geq 0,\varepsilon y_{N}\leq S_{m}\leq
y_{N}\right) &=&\int_{\varepsilon y_{N}}^{y_{N}}\mathbb{P}_{w}\left(
L_{m}\geq 0,S_{m}\in dq\right) \\
&\sim &\int_{\varepsilon y_{N}}^{y_{N}}\frac{1}{a_{m}}g\left( \frac{q-w}{%
a_{m}}\right) C\left( \frac{w}{a_{m}},\frac{q}{a_{m}}\right) dq
\end{eqnarray*}%
uniformly in $w\in \lbrack t_{0}a_{m},t_{1}a_{m}]$. Further, by (\ref{Don00}%
)
\begin{equation}
\frac{\mathbb{P}_{x_{N}}\left( S_{N-m}\in dw,L_{N-m}\geq 0\right) }{dw}\sim
\frac{g(0)}{Na_{N}}V^{-}(x_{N})V^{+}(w)  \label{Addit1}
\end{equation}%
uniformly in $x_{N}=o(a_{N-m})=o(a_{N}),w=o(a_{N-m})$. Since $m=o(N)$, we
have%
\begin{eqnarray*}
&&\mathbb{P}_{x_{N}}\left( \varepsilon y_{N}\leq S_{N-m}\leq
za_{m},L_{N}\geq 0,\varepsilon y_{N}\leq S_{N}\leq y_{N}\right) \\
&&\quad\sim\frac{g(0)V^{-}(x_{N})}{Na_{N}}\int_{\varepsilon
y_{N}}^{za_{m}}V^{+}(w)dw\int_{\varepsilon y_{N}}^{y_{N}}\frac{1}{a_{m}}%
g\left( \frac{q-w}{a_{m}}\right) C\left( \frac{w}{a_{m}},\frac{q}{a_{m}}%
\right) dq \\
&&\quad\sim\frac{g(0)V^{-}(x_{N})a_{m}}{Na_{N}}\int_{\varepsilon
T}^{z}V^{+}(sa_{m})ds\int_{\varepsilon T}^{T}g\left( r-s\right) C\left(
s,r\right) dr \\
&&\quad=\frac{g(0)V^{-}(x_{N})a_{m}V^{+}(a_{m})}{Na_{N}}\int_{\varepsilon
T}^{z}\frac{V^{+}(sa_{m})}{V^{+}(a_{m})}ds\int_{\varepsilon T}^{T}g\left(
r-s\right) C\left( s,r\right) dr \\
&&\quad\sim\frac{g(0)V^{-}(x_{N})a_{m}V^{+}(a_{m})}{Na_{N}}\int_{\varepsilon
T}^{z}s^{\alpha \rho}ds\int_{\varepsilon T}^{T}g(r-s)C\left( s,r\right)dr.
\end{eqnarray*}%
Using (\ref{EquivT}), the equivalence
\begin{equation*}
\left( \alpha \rho +1\right) \int_{0}^{y_{N}}V^{+}(w)dw\sim
y_{N}V^{+}(y_{N}),
\end{equation*}%
the representation%
\begin{equation*}
\mathbb{P}_{x_{N}}\left( L_{N}\geq 0,S_{N}\leq y_{N}\right) =(1+o(1))\frac{%
g(0)}{Na_{N}}V^{-}(x_{N})\int_{0}^{y_{N}}V^{+}(w)dw,
\end{equation*}%
which follows from (\ref{Addit1}), and letting $\varepsilon \downarrow 0$ we
see that, as $N\rightarrow \infty $%
\begin{equation*}
\frac{\mathbb{P}_{x_{N}}\left( S_{N-m}\leq za_{m},L_{N}\geq 0,S_{N}\leq
y_{N}\right) }{\mathbb{P}_{x_{N}}\left( L_{N}\geq 0,S_{N}\leq y_{N}\right) }%
\sim \frac{\alpha \rho +1}{T^{\alpha \rho +1}}\int_{0}^{z}s^{\alpha \rho
}ds\int_{0}^{T}g\left( r-s\right) C\left( s,r\right) dr
\end{equation*}%
uniformly in $x_{N}=o(a_{N})$ and $y_{N}\sim Ta_{m}=o(a_{N})$.

Relation (\ref{Con333}) is proved.

\section{Limit theorem for sequences convergent almost surely\label{s5}}

In this section we need one more conditional law, $\mathbb{P}_{x}^{\uparrow
}\left( \cdot \right) $, the law of the random walk under $\mathbb{P}$
started at $x$ and conditioned to stay non-negative for all times (see, for
example, \cite{BD94}, \cite{CC2008} and \cite{VW09}). It is specified for $%
x\geq 0$ by setting for all $N\in \mathbb{N}$ and any set $B$ belonging to
the $\sigma $-algebra generated by the random variables $S_{1},...,S_{N}$ by
the relation
\begin{equation}
\mathbb{P}_{x}^{\uparrow }(B):=\frac{1}{V^{-}(x)}\mathbb{E}_{x}\left[
V^{-}(S_{N})I\left\{ B\right\} ;L_{N}\geq 0\right] .  \label{StayNonneg}
\end{equation}

The next theorem is a generalization of Lemma 2.5 in \cite{agkv} and Lemma 4
in \cite{VD2022}, where the measure
\begin{equation*}
\mathbb{P}_{x}^{+}(B):=\frac{1}{\hat{V}^{-}(x)}\mathbb{E}_{x}\left[ \hat{V}%
^{-}(S_{N})I\left\{ B\right\} ;L_{N}\geq 0\right]
\end{equation*}%
was used in the statements, which, in view of (\ref{RenewRelation}),
coincides with $\mathbb{P}_{x}^{\uparrow }$.

We take a function $\varphi (n)\rightarrow \infty $ as $n\rightarrow \infty $
and intrudice conditional expactations%
\begin{equation*}
I_{n}(z,m,\varphi ):=\mathbb{E}\left[ H_{n};S_{n-m}\leq za_{m}|S_{n}\leq
\varphi (n),L_{n}\geq 0\right] ,z\in (0,\infty ),
\end{equation*}%
and%
\begin{equation*}
I_{n}^{\ast }(z,m,\varphi ):=\mathbb{E}\left[ H_{n};S_{n-m}-S_{n}\leq
za_{m}|S_{n}\leq \varphi (n),L_{n}\geq 0\right] ,z\in (-\infty ,\infty ).
\end{equation*}

\begin{theorem}
\label{T_cond} Assume Hypothesis $A1$ and A2. Let $H_{1},H_{2},...,$ be a
uniformly bounded sequence of random variables adapted to the filtration $%
\mathcal{\tilde{F}=}\left\{ \mathcal{\tilde{F}}_{k},k\in \mathbb{N}\right\} $
and converging $\mathbb{P}_{0}^{\uparrow }$-a.s. to a random variable $%
H_{\infty }$ as $n\rightarrow \infty $. Suppose that the parameter $%
m=m(n)\rightarrow \infty $ as $n\rightarrow \infty $ in such a way that $%
m=o(n)$. Then

1) if $\varphi (n)=o(a_{m}),$ then
\begin{equation}
\lim_{n\rightarrow \infty }I_{n}(z,m,\varphi )=A_{1}(z)\mathbb{E}%
_{0}^{\uparrow }\left[ H_{\infty }\right] ;  \label{Cond1}
\end{equation}%
2) if $\varphi (n)\sim Ta_{m},T\in (0,\infty ),$ then
\begin{equation}
\lim_{n\rightarrow \infty }I_{n}(z,m,\varphi )=B(z,T)\mathbb{E}%
_{0}^{\uparrow }\left[ H_{\infty }\right] ;  \label{Cond2}
\end{equation}%
3) if $m\rightarrow \infty $ and $a_{m}=o(\varphi (n)),$ then
\begin{equation}
\lim_{n\rightarrow \infty }I_{n}^{\ast }(z,m,\varphi )=\mathbf{P}\left(
Y_{1}\leq z\right) \mathbb{E}_{0}^{\uparrow }\left[ H_{\infty }\right] .
\label{Cond3}
\end{equation}
\end{theorem}

\textbf{Proof}. We prove \ (\ref{Cond1}) for the $(1;0)$-lattice case. For
fixed $1\leq k<n$ and $z\in (0,\infty )$ we consider the quantity%
\begin{eqnarray*}
&&\mathbb{E}\left[ H_{k};S_{n-m}\leq za_{m}|S_{n}\leq \varphi (n),L_{n}\geq 0%
\right] \\
&&\quad =\mathbb{E}\left[ H_{k}\frac{\mathbb{P}_{S_{k}}(S_{n-m-k}^{\prime
}\leq za_{m},S_{n-k}^{\prime }\leq \varphi (n),L_{n-k}^{\prime }\geq 0)}{%
\mathbb{P}\left( S_{n}\leq \varphi (n),L_{n}\geq 0\right) };L_{k}\geq 0%
\right] ,
\end{eqnarray*}%
where $\mathcal{S}^{\prime }=\left\{ S_{n}^{\prime },n=0,1,2,...\right\} $
is a probabilistic copy of the random walk $\mathcal{S}$ being independent
of the set$\mathcal{\ }\left\{ S_{j},j=0,1,...,k\right\} $. We know by \
Corollaries \ref{C_Ysmall} and \ref{C_Yequal} and (\ref{Rough1}) that there
exist constants $C,C_{1},C_{2}$ and $C_{3}$ such that for any fixed $k$ and
all $n\geq k$ and $z>0$
\begin{eqnarray*}
&&\frac{\mathbb{P}_{x}(S_{n-m-k}^{\prime }\leq za_{m},S_{n-k}^{\prime }\leq
\varphi (n),L_{n-k}^{\prime }\geq 0)}{\mathbb{P}\left( S_{n}\leq \varphi
(n),L_{n}\geq 0\right) } \\
&&\qquad \qquad \leq \frac{\mathbb{P}_{x}(S_{n-k}^{\prime }\leq \varphi
(n),L_{n-k}^{\prime }\geq 0)}{\mathbb{P}\left( S_{n}\leq \varphi
(n),L_{n}\geq 0\right) } \\
&&\qquad \qquad \leq \frac{C_{1}\,b_{n-k}\,V^{-}(x)\sum_{j=0}^{\varphi (n)}%
\text{$V$}^{+}(j)}{Cb_{n}\sum_{j=0}^{\varphi (n)}\text{$V$}^{+}(j)}\leq
C_{3}V^{-}(x).
\end{eqnarray*}%
Further, recalling Corollary\textbf{\ }\ref{C_Ysmall}, relation (\ref%
{BasicAsympt3})\ and definition (\ref{Def-b}), we see that, for each fixed $%
x\geq 0$ and $k\in \mathbb{N}$
\begin{equation}
\lim_{n\rightarrow \infty }\frac{\mathbb{P}_{x}(S_{n-m-k}^{\prime }\leq
za_{m},S_{n-k}^{\prime }\leq \varphi (n),L_{n-k}^{\prime }\geq 0)}{\mathbb{P}%
\left( S_{n}\leq \varphi (n),L_{n}\geq 0\right) }=(1-\zeta )A_{1}(z)V^{-}(x).%
\text{ }  \label{Shift}
\end{equation}%
We know by (\ref{V_zero}) that
\begin{eqnarray*}
\mathbb{E}\left[ H_{k}V^{-}(S_{k});L_{k}\geq 0\right] &=&\frac{1}{1-\zeta }%
\times \frac{1}{V^{-}(0)}\mathbb{E}\left[ H_{k}V^{-}(S_{k});L_{k}\geq 0%
\right] \\
&\mathbf{=}&\frac{1}{1-\zeta }\mathbb{E}_{0}^{\uparrow }\left[ H_{k}\right]
<\infty .
\end{eqnarray*}%
Thus, we may apply the dominated convergence theorem and (\ref{Shift}) to
conclude that, for each fixed $k$%
\begin{eqnarray*}
&&\lim_{n\rightarrow \infty }\mathbb{E}\left[ H_{k}\frac{\mathbb{P}%
_{S_{k}}(S_{n-m-k}^{\prime }\leq za_{m},S_{n-k}^{\prime }\leq \varphi
(n),L_{n-k}^{\prime }\geq 0)}{\mathbb{P}\left( S_{n-m}\leq za_{m},S_{n}\leq
\varphi (n),L_{n}\geq 0\right) };L_{k}\geq 0\right] \\
&&\quad =\mathbb{E}\left[ H_{k}\times \lim_{n\rightarrow \infty }\frac{%
\mathbb{P}_{S_{k}}(S_{n-m-k}^{\prime }\leq za_{m},S_{n-k}^{\prime }\leq
\varphi (n),L_{n-k}^{\prime }\geq 0)}{\mathbb{P}\left( S_{n-m}\leq
za_{m},S_{n}\leq \varphi (n),L_{n}\geq 0\right) };L_{k}\geq 0\right] \\
&&\quad =A_{1}(z)\frac{1}{V^{-}(0)}\mathbb{E}\left[ H_{k}V^{-}(S_{k});L_{k}%
\geq 0\right] =A_{1}(z)\mathbb{E}_{0}^{\uparrow }\left[ H_{k}\right] .
\end{eqnarray*}%
To avoid combersome formulas we introduce for $\lambda \geq 1$ a temporary
notation%
\begin{equation*}
\Psi \left( \lambda n,\lambda m,z,\varphi \right) :=\left\{ S_{\lambda
(n-m)}\leq z,S_{\lambda n}\leq \varphi (n),L_{n\lambda }\geq 0\right\} .
\end{equation*}%
In view of (\ref{Rough1}) we have for each $\lambda >1$:
\begin{eqnarray*}
&&\left\vert \mathbb{E}\left[ \left( H_{n}-H_{k}\right) ;\Psi \left( \lambda
n,\lambda m,za_{m},\varphi \right) \right] \right\vert \\
&&\qquad \leq \mathbb{E}\left[ \left\vert H_{n}-H_{k}\right\vert ;S_{\lambda
n}\leq \varphi (n),L_{n\lambda }\geq 0\right] \\
&&\qquad =\mathbb{E}\left[ \left\vert H_{n}-H_{k}\right\vert \mathbb{P}%
_{S_{n}}(S_{(\lambda -1)n}^{\prime }\leq \varphi (n),L_{n(\lambda
-1)}^{\prime }\geq 0);L_{n}\geq 0\right] \\
&&\qquad \leq Cb_{n(\lambda -1)}\sum_{z=0}^{\varphi (n)}V^{+}(z)\times
\mathbb{E}\left[ \left\vert H_{n}-H_{k}\right\vert
\,\,V^{-}(S_{n}),L_{n}\geq 0\right] \\
&&\qquad =Cb_{n(\lambda -1)}\sum_{z=0}^{\varphi (n)}V^{+}(z)\times \frac{1}{%
1-\zeta }\mathbb{E}_{0}^{\uparrow }\left[ \left\vert H_{n}-H_{k}\right\vert
\,\right] .
\end{eqnarray*}%
Further, by Corollary \ref{C_Ysmall}\ and the equivalence $a_{\lambda m}\sim
\lambda ^{1/\alpha }a_{m}$ valid as $m\rightarrow \infty $, we conclude
that, as $m,n\rightarrow \infty $ and $m=o(n)$%
\begin{equation*}
\mathbb{P}\left( \Psi \left( \lambda n,\lambda m,za_{m},\varphi \right)
\right) \sim A_{1}(z\lambda ^{-1/\alpha })\mathbb{P}\left( S_{n\lambda }\leq
\varphi (n),L_{n\lambda }\geq 0\right) .
\end{equation*}%
Hence, using (\ref{Def-b}) and (\ref{BasicAsympt2}) with $x_{N}=0$ we
conclude that%
\begin{eqnarray*}
&&\frac{\left\vert \mathbb{E}\left[ (H_{n}-H_{k});\Psi \left( \lambda
n,\lambda m,za_{m},\varphi \right) \right] \right\vert }{\mathbb{P}\left(
\Psi \left( \lambda n,\lambda m,za_{m},\varphi \right) \right) } \\
&&\qquad \qquad \leq C\mathbb{E}^{\uparrow }\left[ \left\vert
H_{n}-H_{k}\right\vert \right] \frac{b_{n(\lambda -1)}\sum_{z=0}^{\varphi
(n)}V^{+}(z)}{C_{1}\,b_{n\lambda }\,\sum_{z=0}^{\varphi (n)}V^{+}(z)} \\
&&\qquad \qquad \leq C_{2}\left( \frac{\lambda }{\lambda -1}\right)
^{1+1/\alpha }\mathbb{E}_{0}^{\uparrow }\left[ \left\vert
H_{n}-H_{k}\right\vert \right] .
\end{eqnarray*}%
Letting first $n$ and then $k$ to infinity we see that, for each $\lambda >1$
the right-hand side of the previous relation vanishes by the dominated
convergence theorem.

Applying this result we see that%
\begin{eqnarray*}
&&\lim_{n\rightarrow \infty }\mathbb{E}\left[ H_{n}|\Psi \left( \lambda
n,\lambda m,za_{m},\varphi \right) \right] \\
&&\qquad =\lim_{k\rightarrow \infty }\lim_{n\rightarrow \infty }\frac{%
\mathbb{E}\left[ \left( H_{n}-H_{k}\right) ;\Psi \left( \lambda n,\lambda
m,za_{m},\varphi \right) \right] }{\mathbb{P}\left( \Psi \left( \lambda
n,\lambda m,za_{m},\varphi \right) \right) } \\
&&\qquad \quad +\lim_{k\rightarrow \infty }\lim_{n\rightarrow \infty }\frac{%
\mathbb{E}\left[ H_{k};\Psi \left( \lambda n,\lambda m,za_{m},\varphi
\right) \right] }{\mathbb{P}\left( \Psi \left( \lambda n,\lambda
m,za_{m},\varphi \right) \right) } \\
&&\qquad =\lim_{k\rightarrow \infty }A_{1}(z\lambda ^{-1/\alpha })\mathbb{E}%
_{0}^{\uparrow }\left[ H_{k}\right] =A_{1}(z\lambda ^{-1/\alpha })\mathbb{E}%
_{0}^{\uparrow }\left[ H_{\infty }\right] ,
\end{eqnarray*}%
which we rewrite as
\begin{eqnarray*}
&&\mathbb{E}\left[ H_{n};\Psi \left( \lambda n,\lambda m,za_{m},\varphi
\right) \right] \\
&=&\left( A_{1}(z\lambda ^{-1/\alpha })\mathbb{E}_{0}^{\uparrow }\left[
H_{\infty }\right] +o(1)\right) \mathbb{P}\left( \Psi \left( \lambda
n,\lambda m,za_{m},\varphi \right) \right) .
\end{eqnarray*}%
Assuming without loss of generality that $H_{\infty }>0$ \ and \ $%
A_{1}(z\lambda ^{-1/\alpha })\mathbb{E}_{0}^{\uparrow }\left[ H_{\infty }%
\right] \leq 1,$ we conclude that
\begin{eqnarray*}
&&|\mathbb{E}\left[ H_{n};\Psi \left( n,m,za_{m},\varphi \right) \right]
-A_{1}(z\lambda ^{-1/\alpha })\mathbb{E}_{0}^{\uparrow }\left[ H_{\infty }%
\right] \mathbb{P}\left( \Psi \left( \lambda n,\lambda m,za_{m},\varphi
\right) \right) | \\
&&\quad \leq |\mathbb{E}\left[ H_{n};\Psi \left( \lambda n,\lambda
m,za_{m},\varphi \right) \right] -A_{1}(z\lambda ^{-1/\alpha })\mathbb{E}%
_{0}^{\uparrow }\left[ H_{\infty }\right] \mathbb{P}\left( \Psi \left(
\lambda n,\lambda m,za_{m},\varphi \right) \right) | \\
&&\qquad +|\mathbb{P}\left( \Psi \left( \lambda n,\lambda m,za_{m},\varphi
\right) \right) -\mathbb{P}\left( \Psi \left( n,m,za_{m},\varphi \right)
\right) |.
\end{eqnarray*}%
We have proved that the first difference at the right-hand side of the
inequality is of the order

\begin{equation*}
o\left( \mathbb{P}\left( \Psi \left( \lambda n,\lambda m,za_{m},\varphi
\right) \right) \right)
\end{equation*}%
as $n\rightarrow \infty $, and, therefore, of the order $o\left( \mathbb{P}%
\left( S_{n}\leq \varphi (n),L_{n}\geq 0\right) \right) ,$ since
\begin{equation}
\lim_{n\rightarrow \infty }\frac{\mathbb{P}\left( S_{n\lambda }\leq \varphi
(n),L_{n\lambda }\geq 0\right) }{\mathbb{P}\left( S_{n}\leq \varphi
(n),L_{n}\geq 0\right) }=\lim_{n\rightarrow \infty }\frac{b_{n\lambda }}{%
b_{n}}=\lambda ^{1+1/\alpha }  \label{RatioB}
\end{equation}%
by (\ref{BasicAsympt2}) and (\ref{Def-b}).

Further, again by (\ref{BasicAsympt2}), Corollary \ref{C_Ysmall}, and
definition (\ref{Def-b}) we have
\begin{eqnarray*}
&&|\mathbb{P}\left( \Psi \left( \lambda n,\lambda m,za_{m},\varphi \right)
\right) -\mathbb{P}\left( \Psi \left( n,m,za_{m},\varphi \right) \right) | \\
&\leq &\left\vert \mathbb{P}\left( \Psi \left( \lambda n,\lambda
m,za_{m},\varphi \right) \right) -g(0)A_{1}(z\lambda ^{-1/\alpha
})b_{n\lambda }\sum_{j=0}^{\varphi (n)}V^{+}(j)\right\vert \\
&&+\left\vert \mathbb{P}\left( \Psi \left( n,m,za_{m},\varphi \right)
\right) -g(0)A_{1}(z)b_{n}\sum_{j=0}^{\varphi (n)}V^{+}(j)\right\vert \\
&&+g(0)\left\vert A_{1}(z\lambda ^{-1/\alpha })b_{n\lambda
}-A_{1}(z)b_{n}\right\vert \sum_{j=0}^{\varphi (n)}V^{+}(j) \\
&=&o\left( b_{n}\sum_{j=0}^{\varphi (n)}V^{+}(j)\right)
+g(0)b_{n}A_{1}(z)\left\vert \frac{A_{1}(z\lambda ^{-1/\alpha })b_{n\lambda }%
}{A_{1}(z)b_{n}}-1\right\vert \sum_{j=0}^{\varphi (n)}V^{+}(j).
\end{eqnarray*}%
Hence, letting $\lambda \downarrow 1$, using (\ref{RatioB}) and the
continuity of $A_{1}(z)$ we conclude that%
\begin{equation*}
\lim_{\lambda \downarrow 1}\lim_{n\rightarrow \infty }\frac{\left\vert
\mathbb{P}\left( \Psi \left( \lambda n,\lambda m,za_{m},\varphi \right)
\right) -\mathbb{P}\left( \Psi \left( n,m,za_{m},\varphi \right) \right)
\right\vert }{b_{n}\sum_{j=0}^{\varphi (n)}V^{+}(j)}=0.
\end{equation*}%
Combining the obtained estimates we get (\ref{Cond1}) for the $(1;0)$%
-lattice case.

To prove (\ref{Cond2}) for the $(1;0)$-lattice case one may repeat almost
literally the argument used to justify \bigskip (\ref{Cond1}) replacing the
reference to Corollary \ref{C_Ysmall} by the reference to Corollary \ref%
{C_Yequal}.

To prove (\ref{Cond3}) for the $(1;0)$-lattice case one can repeat almost
literally the argument used to justify \bigskip (\ref{Cond1}) replacing
there the reference to Corollary \ref{C_Ysmall} by the reference to
Corollary \ref{C_Ybig}.

To prove (\ref{Cond1}) for the absolutely continuous case it is necessary to
replace in the arguments above $\sum_{j=0}^{\varphi (n)}V^{+}(j)$ by $%
\int_{0}^{\varphi (n)}V^{+}(w)dw$ and to use Lemma \ref{L_TotalCont}.
Similar arguments allow to establish (\ref{Cond2}) and (\ref{Cond3}) for the
absolutely continuous distribution of $X_{1}$.

\section{Branching processes in random environment\label{s6}}

In this section we apply the results obtained for random walks to study the
population size of the critical branching processes evolving in unfavorable
random environments. To describe the model and the problems we plan to
consider we intruduce the space $\mathcal{M}$ $=\left\{ \mathfrak{f}\right\}
$ of all probability measures on $\mathbb{N}_{0}$. For notational reason, we
identify a measure $\mathfrak{f}=\left\{ \mathfrak{f}(\left\{ 0\right\} ),%
\mathfrak{f}(\left\{ 1\right\} ),...\right\} \in $ $\mathcal{M}$ with the
respective probability generating function%
\begin{equation*}
f(s)=\sum_{k=0}^{\infty }\mathfrak{f}(\left\{ k\right\} )s^{k},\quad s\in
\lbrack 0,1],
\end{equation*}%
and make no difference between $\mathfrak{f}$ and $f$. Equipped with a
metric, $\mathcal{M}$ $=\left\{ \mathfrak{f}\right\} =\left\{ f\right\} $
becomes a Polish space. Let
\begin{equation*}
F(s)=\sum_{j=0}^{\infty }F\left( \left\{ j\right\} \right) s^{j},\quad s\in
\lbrack 0,1],
\end{equation*}%
be a random variable taking values in $\mathcal{M}$, and let
\begin{equation*}
F_{n}(s)=\sum_{j=0}^{\infty }F_{n}\left( \left\{ j\right\} \right)
s^{j},\quad s\in \lbrack 0,1],\quad n\in \mathbb{N}:=\mathbb{N}%
_{0}\backslash \left\{ 0\right\} ,
\end{equation*}%
be a sequence of independent copies of $F$. The infinite sequence $\mathcal{E%
}=\left\{ F_{n},n\in \mathbb{N}\right\} $ is called a random environment.

A sequence of nonnegative random variables $\mathcal{Z}=\left\{ Z_{n},\ n\in
\mathbb{N}_{0}\right\} $ specified on\ a probability space $(\Omega ,%
\mathcal{F},\mathbb{P})$ is called a branching process in random environment
(BPRE), if $Z_{0}$ is independent of $\mathcal{E}$ and, given $\mathcal{E}$
the process $\mathcal{Z}$ is a Markov chain with
\begin{equation*}
\mathcal{L}\left( Z_{n}|Z_{n-1}=z_{n-1},\mathcal{E}=(f_{1},f_{2},...)\right)
=\mathcal{L}(\xi _{n1}+\ldots +\xi _{ny_{n-1}})
\end{equation*}%
for all $n\in \mathbb{N}$, $z_{n-1}\in \mathbb{N}_{0}$ and $%
f_{1},f_{2},...\in \mathcal{M}$, where $\xi _{n1},\xi _{n2},\ldots $ is a
sequence of i.i.d. random variables with distribution $f_{n}.$ Thus, $%
Z_{n-1} $ is the $(n-1)$th generation size of the population of the
branching process and $f_{n}$ is the offspring distribution of an individual
at generation $n-1$.

The sequence
\begin{equation*}
S_{0}=0,\quad S_{n}=X_{1}+...+X_{n},\ n\geq 1,
\end{equation*}%
where $X_{i}=\log F_{i}^{\prime }(1),i=1,2,...$ is called the associated
random walk for the process $\mathcal{Z}$.

We impose the following restrictions on the properties of the BPRE.

\paragraph{Condition B1.}

\emph{The elements of the associated random walk satisfy Hypothesis A1 and
A2. }

According to the classification of \ BPREs (see, for instance, \cite{agkv}
and \cite{KV2017}), Condition B1 means that we consider the critical BPRE's.

Our second assumption on the environment concerns reproduction laws of
particles. Set%
\begin{equation*}
\gamma (b)=\frac{\sum_{k=b}^{\infty }k^{2}F\left( \left\{ k\right\} \right)
}{\left( \sum_{i=b}^{\infty }iF\left( \left\{ i\right\} \right) \right) ^{2}}%
.
\end{equation*}

\paragraph{Condition B2.}

\emph{There exist $\varepsilon >0$ and $b\in $}$\mathbb{N}$ \emph{such that\
} \emph{\ }
\begin{equation*}
\mathbb{E}[(\log ^{+}\gamma (b))^{\alpha +\varepsilon }]\ <\ \infty \ ,
\end{equation*}%
\emph{where }$\log ^{+}x=\log (\max (x,1))$\emph{.}

It is known (see \cite[Theorem 1.1 and Corollary 1.2]{agkv}) that if
Conditions B1, B2 are valid then there exist a number $\theta \in (0,\infty
) $ and a sequence $l(1),l(2)...,$ slowly varying at infinity such that, as $%
n\rightarrow \infty $%
\begin{equation*}
\mathbb{P}\left( Z_{n}>0\right) \sim \theta n^{-(1-\rho )}l(n)
\end{equation*}%
and for any $t\in \left[ 0,1\right] $ and any $x\geq 0$
\begin{equation}
\lim_{n\rightarrow \infty }\mathbb{P}\left( \frac{\log Z_{\left[ nt\right] }%
}{a_{n}}\leq x|Z_{n}>0\right) =\lim_{n\rightarrow \infty }\mathbb{P}\left(
\frac{S_{\left[ nt\right] }}{a_{n}}\leq x|Z_{n}>0\right) =\mathbf{P}\left(
Y_{t}^{+}\leq x\right) ,  \label{Meander0}
\end{equation}%
where $\mathcal{Y}^{+}=\left\{ Y_{t}^{+},0\leq t\leq 1\right\} $ denotes the
meander of the strictly stable process $\mathcal{Y}$ with index $\alpha $.

Thus, if a BPRE is critical, then, given $Z_{n}>0$ the random variables $%
\log Z_{\left[ nt\right] },t\in (0,1],$ and $S_{n},$ the value of the
associated random walk that provides survival of the population to a distant
moment $n$, \ grow like $a_{n}$ times random positive multipliers.

These results were complemented in \cite{VD2023} by considering the properly
scaled distributions of the random variables $\log Z_{\left[ nt\right]
},t\in (0,1]$ given that $Z_{n}>0$ and $S_{n}\leq \varphi (n),$ where $%
\varphi (n)\rightarrow \infty $ as $n\rightarrow \infty $ in such a way that
$\varphi (n)=o(a_{n})$. In view of (\ref{Meander0}) the event $\left\{
S_{n}\leq \varphi (n)\right\} $ may be treated in this case as an
unfavorable one for the development of a critical branching process given
its survival.

Introduce the notation $A_{u.s}:=\left\{ Z_{n}>0\text{ for all }n>0\right\} $
for the event of ultimate survival. Observe that according to Theorem 1 in
\cite{VD2022},\ if Conditions B1, B2 are valid, the distribution of the
increments of the associated random walk is absolutely continuous and $%
\varphi (n)=o(a_{n})$, then
\begin{eqnarray}
\mathbb{P}\left( S_{n}\leq \varphi (n),Z_{n}>0\right) &\sim &\Theta \mathbb{P%
}\left( S_{n}\leq \varphi (n),L_{n}>0\right)  \notag \\
&\sim &\frac{\Theta }{na_{n}}\int_{0}^{\varphi (n)}V^{+}(w)dw,
\label{AsymMain1}
\end{eqnarray}%
as $n\rightarrow \infty $, where%
\begin{equation}
\Theta =\sum_{j=0}^{\infty }\sum_{k=1}^{\infty }\mathbb{P}(Z_{j}=k,\tau
_{j}=j)\mathbb{P}^{\uparrow }\left( A_{u.s}|Z_{0}=k\right) \in (0,\infty ).
\label{DefTheta}
\end{equation}

The same asymptotics is valid for $(1;0)$-lattice case, where $%
\int_{0}^{\varphi (n)}V^{+}(w)dw$ should be replaced by $\sum_{j=0}^{\varphi
(n)}V^{+}(j).$

Along with the asymptotic behavior of the survival probability of a critical
BPRE, the growth of the population size of such a process at the logarithmic
scale was described in the following theorem.

\begin{theorem}
\label{T_smallDevi}(see \cite{VD2023}) Let Conditions B1, B2 be valid. If $%
\varphi (n)\rightarrow \infty $ as $n\rightarrow \infty $ in such a way that
$\varphi (n)=o(a_{n})$ then for any $y\in (0,1]$%
\begin{equation}
\lim_{n\rightarrow \infty }\mathbb{P}\left( \frac{1}{\varphi (n)}\log
Z_{n}\leq y|S_{n}\leq \varphi (n),Z_{n}>0\right) =y^{\alpha \rho +1}
\label{For1}
\end{equation}%
and for any $t\in (0,1)$ and $x\in \lbrack 0,\infty )$%
\begin{equation}
\lim_{n\rightarrow \infty }\mathbb{P}\left( \frac{1}{a_{n}}\log Z_{\left[ nt%
\right] }\leq x|S_{n}\leq \varphi (n),Z_{n}>0\right) =\mathbf{P}\left(
Y_{t}^{++}\leq x\right) ,  \label{ForLess1}
\end{equation}%
where $\mathcal{Y}^{++}=\left\{ Y_{t}^{++},0\leq t\leq 1\right\} $ denotes
the excursion of the strictly stable process $\mathcal{Y}$ with index $%
\alpha $.
\end{theorem}

\bigskip There is a striking difference between the orders of scaling in (%
\ref{For1}) and (\ref{ForLess1}) indicating that there should be a phase
transition in the growth rate of the population size when we consider the
process within the interval $[n-m,n)$, $m=o(n),m\rightarrow \infty $. This
is indeed the case and it is the aim of this section to demonstrate that if
we scale $\log Z_{n-m}$ by $a_{m}$, then three different limiting
distributions appear as $\min (n-m,m)\rightarrow \infty $. The form of these
distributions depends on which of the following three asymptotic relations
is fulfilled: $\varphi (n)=o(m),\varphi (n)$ is proportional to $m$, and $%
m=o(\varphi (n))$.

Our main result looks as follows. Recall that in the lattice case we
restrict ourselves by the distributions concentrated in the $(1;0)$-lattice.

\begin{theorem}
\label{T_global}Let Conditions B1, B2 be valid, $\min (m,n)\rightarrow
\infty $ and $m=o(n)$. Then

1) If $\varphi (n)\rightarrow \infty $ as $n\rightarrow \infty $ in such a
way that $\varphi (n)=o(a_{m}),$ then for any $z\in (0,\infty )$
\begin{equation*}
\lim_{n\rightarrow \infty }\mathbb{P}\left( \frac{1}{a_{m}}\log Z_{n-m}\leq
z|S_{n}\leq \varphi (n),Z_{n}>0\right) =A_{1}(z),
\end{equation*}%
where $A_{1}(z)$ is defined in (\ref{DefA1(z)});

2) If $\varphi (n)\rightarrow \infty $ as $n\rightarrow \infty $ in such a
way that $\varphi (n)\sim Ta_{m},T\in (0,\infty ),$ then for any $z\in
\lbrack 0,\infty )$
\begin{equation*}
\lim_{n\rightarrow \infty }\mathbb{P}\left( \frac{1}{a_{m}}\log Z_{n-m}\leq
z|S_{n}\leq \varphi (n),Z_{n}>0\right) =B(z,T),
\end{equation*}%
where $B(z,T)$ is defined in (\ref{DefB(z,T)});

3) if $m=o(\varphi (n)),\varphi (n)=o(a_{n}),$ then for any $z\in \lbrack
0,\infty )$
\begin{equation*}
\lim_{n\rightarrow \infty }\mathbb{P}\left( \frac{1}{a_{m}}\left( \log
Z_{n-m}-S_{n}\right) \leq z|S_{n}\leq \varphi (n),Z_{n}>0\right) =\mathbf{P}%
\left( Y_{1}\leq z\right) .
\end{equation*}
\end{theorem}

\textbf{\ \ Proof. }For integers $0\leq r\leq n$, consider the rescaled
generation size process $\mathcal{X}^{r,n}=\left\{ X_{t}^{r,n},0\leq t\leq
1\right\} $, given by%
\begin{equation*}
X_{t}^{r,n}=e^{-S_{r+\left[ (n-r)t\right] }}Z_{r+\left[ (n-r)t\right]
},0\leq t\leq 1.\quad
\end{equation*}%
It was shown in Theorem 1 of \cite{VD2023} that if $r_{1},r_{2},...,$ is a
sequence of positive integers such that $r_{n}\leq n$ and $r_{n}\rightarrow
\infty $ and $\varphi (n)\rightarrow \infty $ as $n\rightarrow \infty $ in
such a way that $\varphi (n)=o(a_{n})$, then, as $n\rightarrow \infty $
\begin{equation}
\mathcal{L}\left\{ X_{t}^{r_{n},n},0\leq t\leq 1|S_{n}\leq \varphi
(n),Z_{n}>0\right\} \Longrightarrow \mathcal{L}\left\{ W_{t},0\leq t\leq
1\right\} ,  \label{Skoroh2}
\end{equation}%
where $W_{t},0\leq t\leq 1,$ is a stochastic process with a.s. constant
paths, that is%
\begin{equation*}
\mathbf{P}\left( W_{t}=W\text{ for all }t\in (0,1]\right) =1
\end{equation*}%
for some random variable $W$. Furthermore,%
\begin{equation}
\mathbf{P}\left( 0<W<\infty \right) =1.  \label{RangeW}
\end{equation}%
Put

\begin{equation*}
\hat{Z}(k)=e^{-S_{k}}Z_{k},k=0,1,2,...,
\end{equation*}%
and set for brevity
\begin{equation}
R_{n}(x):=\left\{ S_{n}\leq x,Z_{n}>0\right\} ,\quad Q_{n}(x):=\left\{
S_{n}\leq x,L_{n}\geq 0\right\} .\quad  \label{d2}
\end{equation}%
It follows from (\ref{Skoroh2}) that if $\min (m,n)\rightarrow \infty $ and $%
m=o(n)$, then
\begin{equation*}
\lim_{n\rightarrow \infty }\mathbb{P}\left( \hat{Z}(n-m)\leq x|R_{n}(\varphi
(n))\right) =\mathbf{P}\left( W\leq x\right)
\end{equation*}%
for any point $x\in (0,\infty )$ of continuity of the distribution of the
random variable $W$. Using (\ref{RangeW}) we conclude that, for any $%
\varepsilon >0$ there exists $M=M(m,n)$ such that
\begin{equation}
\mathbb{P}\left( \hat{Z}(n-m)\in (M^{-1},M)|R_{n}(\varphi (n))\right) \geq
1-\varepsilon  \label{Negl11}
\end{equation}%
for all $m\geq m_{0}$ and $n-m\geq n_{0}$. For $z\in (0,\infty )$ write%
\begin{eqnarray}
&&\mathbb{P}\left( \log Z_{n-m}\leq za_{m},R_{n}(\varphi (n))\right)  \notag
\\
&=&\mathbb{P}\left( \log \hat{Z}_{n-m}+S_{n-m}\leq za_{m},\hat{Z}_{n-m}\in
(M^{-1},M),R_{n}(\varphi (n))\right)  \notag \\
&&+\mathbb{P}\left( \log Z_{n-m}\leq za_{m},\hat{Z}_{n}\notin
(M^{-1},M),R_{n}(\varphi (n))\right)  \label{Decompose0}
\end{eqnarray}%
and investigate separately the asymptotic behavior of the summands at the
right-hand side as $\min (m,n)\rightarrow \infty $, $m=o(n)$. In view of (%
\ref{Negl11}) it is sufficient to analyze
\begin{eqnarray*}
&&\lim_{\min (m,n-m)\rightarrow \infty }\frac{\mathbb{P}\left( \log \hat{Z}%
_{n-m}+S_{n-m}\leq za_{m},\hat{Z}_{n-m}\in (M^{-1},M),R_{n}(\varphi
(n))\right) }{\mathbb{P}\left( R_{n}(\varphi (n))\right) } \\
&&\qquad \qquad \qquad \qquad \qquad \qquad \leq \lim_{\min
(m,n-m)\rightarrow \infty }\frac{\mathbb{P}\left( S_{n-m}\leq
za_{m},R_{n}(\varphi (n))\right) }{\mathbb{P}\left( R_{n}(\varphi
(n))\right) }.
\end{eqnarray*}%
We fix $J\in \lbrack 1,n]$ and write%
\begin{eqnarray}
\mathbb{P}\left( S_{n-m}\leq za_{m},R_{n}(\varphi (n))\right) &=&\mathbb{P}%
\left( S_{n-m}\leq za_{m},R_{n}(\varphi (n)),\tau _{n}>J\right)  \notag \\
&&+\mathbb{P}\left( S_{n-m}\leq za_{m},R_{n}(\varphi (n)),\tau _{n}\leq
J\right) .  \label{Decompose1}
\end{eqnarray}%
According to Lemma 5 in \cite{VD2022}
\begin{equation*}
\lim_{J\rightarrow \infty }\limsup_{n\rightarrow \infty }\frac{\mathbb{P}%
\left( R_{n}(\varphi (n)),\tau _{n}>J\right) }{\mathbb{P}\left(
Q_{n}(\varphi (n))\right) }=0.
\end{equation*}%
Using this result and recalling (\ref{AsymMain1}) we obtain
\begin{eqnarray*}
&&\lim_{J\rightarrow \infty }\limsup_{n\rightarrow \infty }\frac{\mathbb{P}%
\left( S_{n-m}\leq za_{m},R_{n}(\varphi (n)),\tau _{n}>J\right) }{\mathbb{P}%
\left( R_{n}(\varphi (n))\right) } \\
&&\qquad \qquad \leq \lim_{J\rightarrow \infty }\limsup_{n\rightarrow \infty
}\frac{\mathbb{P}\left( R_{n}(\varphi (n)),\tau _{n}>J\right) }{\mathbb{P}%
\left( R_{n}(\varphi (n))\right) }=0.
\end{eqnarray*}%
For fixed $j\in \lbrack 1,J]$ we have
\begin{eqnarray*}
&&\mathbb{P}\left( S_{n-m}\leq za_{m},R_{n}(\varphi (n)),\tau
_{n}=j;S_{j}\leq -\sqrt{\varphi (n)}\right) \\
&&\quad \leq \mathbb{P}\left( S_{n}\leq \varphi (n),Z_{j}>0,\tau
_{n}=j;S_{j}\leq -\sqrt{\varphi (n)}\right) \\
&&\quad =\mathbb{E}\left[ \mathbb{P}\left( Z_{j}>0|\mathcal{E}\right)
,S_{n}\leq \varphi (n),\tau _{n}=j;S_{j}\leq -\sqrt{\varphi (n)}\right] \\
&&\quad \leq \mathbb{E}\left[ e^{S_{j}},S_{n}\leq \varphi (n),\tau
_{n}=j;S_{j}\leq -\sqrt{\varphi (n)}\right] \\
&&\qquad =o(\mathbb{P}\left( Q_{n-j}(\varphi (n))\right) ),
\end{eqnarray*}%
where the last equality is established in Lemma 5 in \cite{VD2023}. Further,%
\begin{eqnarray*}
&&\mathbb{P}\left( S_{n-m}\leq za_{m},R_{n}(\varphi (n)),\tau
_{n}=j;Z_{j}>K,S_{j}>-\sqrt{\varphi (n)}\right) \\
&&\quad \leq \mathbb{P}\left( S_{n}\leq \varphi (n),\tau _{n}=j,S_{j}>-\sqrt{%
\varphi (n)},Z_{j}>K\right) \\
&&\quad =\mathbb{P}\left( S_{n}-S_{j}\leq \varphi (n)-S_{j},\tau
_{n}=j,S_{j}>-\sqrt{\varphi (n)},Z_{j}>K\right) \\
&&\quad \leq \mathbb{P}\left( S_{n}-S_{j}\leq \varphi (n)+\sqrt{\varphi (n)}%
,\tau _{n}=j,Z_{j}>K\right) \\
&&\quad =\mathbb{P}\left( \tau _{j}=j;Z_{j}>K\right) \mathbb{P}\left(
S_{n-j}\leq \varphi (n)+\sqrt{\varphi (n)},L_{n-j}\geq 0\right) \\
&&\quad \leq \delta (K)\mathbb{P}\left( S_{n-j}\leq \varphi (n)+\sqrt{%
\varphi (n)},L_{n-j}\geq 0\right) ,
\end{eqnarray*}%
where at the last step we have used the estimates
\begin{equation}
\mathbb{P}\left( \tau _{j}=j,Z_{j}>K\right) \leq \mathbb{P}\left(
Z_{j}>K\right) =\delta (K)\rightarrow 0  \label{LastCommon}
\end{equation}%
as $K\rightarrow \infty $.

Starting from this point and till the end of the proof we assume that the
distribution of $X_{1}$ is absolutely continuous. To prove the desired
statements for the lattice distribution it is necessary everywhere below to
replace everywhere $\int $ by $\sum $.

Consider the right-hand side of the equality
\begin{eqnarray}
&&\mathbb{P}\left( S_{n-m}\leq za_{m},R_{n}(\varphi (n)),\tau _{n}=j,S_{j}>-%
\sqrt{\varphi (n)},Z_{j}\leq K\right)  \notag \\
&&\quad =\int_{-\sqrt{\varphi (n)}}^{0}\sum_{k=1}^{K}\mathbb{P}\left(
S_{j}\in dq,Z_{j}=k,\tau _{j}=j\right)\times  \notag \\
&&\quad \times \mathbb{E}\left[ \mathbb{P}\left( Z_{n-j}>0|\mathcal{E}%
,Z_{0}=k\right) ;S_{n-m-j}\leq za_{m},Q_{n-j}(\varphi (n)-q)\right] .
\label{LargeSj}
\end{eqnarray}

Observe that by monotonicity of the survival probability for each $k\in
\mathbb{N}$
\begin{eqnarray}
H_{n-j}(k) &:&=\mathbb{P}\left( Z_{n-j}>0|\mathcal{E},Z_{0}=k\right)  \notag
\\
&\rightarrow &\mathbb{P}^{\uparrow }\left( A_{u.s}|\mathcal{E}%
,Z_{0}=k\right) =:H_{\infty }\left( k\right)  \label{ConvAs}
\end{eqnarray}%
$\mathbb{P}^{\uparrow }$-a.s. as $n-j\rightarrow \infty $. Moreover, $%
H_{\infty }\left( k\right) >0$ $\mathbb{P}^{\uparrow }$-a.s. according to
Proposition 3.1 in \cite{agkv}. Further, for $q\in (-\sqrt{\varphi (n)},0]$
we have
\begin{eqnarray}
&&\mathbb{E}\left[ H_{n-j}(k);S_{n-m-j}\leq za_{m},Q_{n-j}(\varphi (n)+\sqrt{%
\varphi (n)})\right]  \notag \\
&&\quad \geq \mathbb{E}\left[ H_{n-j}(k);S_{n-m-j}\leq
za_{m},Q_{n-j}(\varphi (n)-q)\right]  \notag \\
&&\quad \geq \mathbb{E}\left[ H_{n-j}(k);S_{n-m-j,n-j}\leq
za_{m},Q_{n-j}(\varphi (n))\right]  \label{Envelope}
\end{eqnarray}%
and, according to Theorem \ref{T_cond}, as $n\rightarrow \infty $
\begin{eqnarray}
&&\mathbb{E}\left[ H_{n-j}(k);S_{n-m-j}\leq za_{m},Q_{n-j}(\varphi (n))%
\right]  \notag \\
&&\qquad \sim \mathbb{P}\left( S_{n-m-j,n-j}\leq za_{m},Q_{n-j}(\varphi
(n))\right) \mathbb{E}^{\uparrow }\left[ H_{\infty }(k)\right]  \notag \\
&&\qquad =\mathbb{P}\left( S_{n-m-j,n-j}\leq za_{m}|Q_{n-j}(\varphi
(n))\right) \times  \notag \\
&&\qquad \quad \times \mathbb{P}\left( Q_{n-j}(\varphi (n))\right) \mathbb{P}%
^{\uparrow }\left( A_{u.s}|Z_{0}=k\right)  \label{EstBelow}
\end{eqnarray}%
and
\begin{eqnarray}
&&\mathbb{E}\left[ H_{n-j}(k);S_{n-m-j}\leq za_{m},Q_{n-j}(\varphi (n)+\sqrt{%
\varphi (n)})\right]  \notag \\
&&\quad \sim \mathbb{P}\left( S_{n-m-j,n-j}\leq za_{m},Q_{n-j}(\varphi (n)+%
\sqrt{\varphi (n)})\right) \mathbb{P}^{\uparrow }\left(
A_{u.s}|Z_{0}=k\right)  \notag \\
&&\quad =\mathbb{P}\left( S_{n-m-j,n-j}\leq za_{m}|Q_{n-j}(\varphi (n)+\sqrt{%
\varphi (n)})\right) \times  \notag \\
&&\quad \quad \times \mathbb{P}\left( Q_{n-j}(\varphi (n)+\sqrt{\varphi (n)}%
)\right) \mathbb{P}^{\uparrow }\left( A_{u.s}|Z_{0}=k\right) .
\label{EstAbove}
\end{eqnarray}

In view of properties of regularly varying functions, the second equivalence
in (\ref{AsymMain1}), and the asymptotic representation%
\begin{equation*}
\int_{0}^{\varphi (n)}V^{+}(w)dw\sim \frac{1}{\alpha \rho +1}\varphi
(n)V^{+}(\varphi (n)),
\end{equation*}%
as $n\rightarrow \infty $, which follows from (\ref{Regular1}), we have
\begin{eqnarray*}
1 &\leq &\lim_{n\rightarrow \infty }\sup_{-\sqrt{\varphi (n)}\leq q\leq 0}%
\frac{\mathbb{P}\left( S_{n-j}\leq \varphi (n)-q,L_{n-j}\geq 0\right) }{%
\mathbb{P}\left( S_{n-j}\leq \varphi (n),L_{n-j}\geq 0\right) } \\
&=&\lim_{n\rightarrow \infty }\frac{\mathbb{P}\left( S_{n-j}\leq \varphi (n)+%
\sqrt{\varphi (n)},L_{n-j}\geq 0\right) }{\mathbb{P}\left( S_{n-j}\leq
\varphi (n),L_{n-j}\geq 0\right) }=1.
\end{eqnarray*}%
This estimate combined with (\ref{d2}) shows that, for all $n\geq j$
\begin{eqnarray*}
&&\sup_{-\sqrt{\varphi (n)}\leq q\leq 0}\frac{\mathbb{E}\left[
H_{n-j}(k);S_{n-m-j}\leq za_{m};Q_{n-j}(\varphi (n)-q)\right] }{\mathbb{P}%
\left( Q_{n-j}(\varphi (n))\right) } \\
&&\qquad \leq \sup_{-\sqrt{\varphi (n)}\leq q\leq 0}\frac{\mathbb{P}\left(
Q_{n-j}(\varphi (n)-q)\right) }{\mathbb{P}\left( Q_{n-j}(\varphi (n))\right)
} \\
&&\qquad \quad =\frac{\mathbb{P}\left( S_{n-j}\leq \varphi (n)+\sqrt{\varphi
(n)},L_{n-j}\geq 0\right) }{\mathbb{P}\left( S_{n-j}\leq \varphi
(n),L_{n-j}\geq 0\right) }\leq C.
\end{eqnarray*}

Now we separately consider the cases $\varphi (n)=o(a_{m})$ and $\varphi
(n)\sim Ta_{m}$ and $a_{m}=o(\varphi (n)).$

1) Assume $\varphi (n)=o(a_{m})$. In this case, setting $H_{n}:=\mathbb{P}%
\left( Z_{n-j}>0|\mathcal{E},Z_{0}=k\right) $ in (\ref{Cond1}) and recalling
(\ref{ConvAs}), we obtain for each $q\in \left[ -\sqrt{\varphi (n)},0\right]
:$
\begin{eqnarray*}
&&\lim_{n\rightarrow \infty }\frac{\mathbb{E}\left[ H_{n-j}(k);S_{n-m-j}\leq
za_{m},Q_{n-j}(\varphi (n)-q)\right] }{\mathbb{P}\left( Q_{n-j}(\varphi
(n))\right) } \\
&&\qquad \qquad \qquad \qquad =A_{1}(z)\mathbb{P}^{\uparrow
}(A_{u.s}|Z_{0}=k)\lim_{n\rightarrow \infty }\frac{\mathbb{P}\left(
Q_{n-j}(\varphi (n)-q)\right) }{\mathbb{P}\left( Q_{n-j}(\varphi (n))\right)
} \\
&&\qquad \qquad \qquad \qquad =A_{1}(z)\mathbb{P}^{\uparrow
}(A_{u.s}|Z_{0}=k).
\end{eqnarray*}%
For $j\in \mathbb{N}_{0}$ and $K\in \mathbb{N}\cup \left\{ \infty \right\} $
set
\begin{equation*}
\Theta _{j}(K):=\sum_{k=1}^{K}\mathbb{P}\left( Z_{j}=k,\tau _{j}=j\right)
\mathbb{P}^{\uparrow }\left( A_{u.s}|Z_{0}=k\right) .
\end{equation*}%
Note that $\Theta _{j}(\infty )\leq \mathbb{P}\left( \tau _{j}=j\right) $.

Applying the dominated convergence theorem we have for fixed $K<\infty $
\begin{eqnarray}
&&\lim_{n\rightarrow \infty }\frac{\mathbb{P}\left( S_{n-m}\leq
za_{m},R_{n}(\varphi (n)),\tau _{n}=j,S_{j}>-\sqrt{\varphi (n)},Z_{j}\leq
K\right) }{\mathbb{P}\left( Q_{n-j}(\varphi (n))\right) }  \notag \\
&&\qquad =\lim_{n\rightarrow \infty }\int_{-\sqrt{\varphi (n)}%
}^{0}\sum_{k=1}^{K}\mathbb{P}\left( S_{j}\in dq,Z_{j}=k,\tau _{j}=j\right)
\notag \\
&&\qquad \qquad \times \frac{\mathbb{E}\left[ H_{n-j}(k);S_{n-m-j}\leq
za_{m},Q_{n-j}(\varphi (n)-q)\right] }{\mathbb{P}\left( Q_{n-j}(\varphi
(n))\right) }  \notag \\
&&\qquad \qquad =A_{1}(z)\int_{-\infty }^{0}\sum_{k=1}^{K}\mathbb{P}\left(
S_{j}\in dq,Z_{j}=k,\tau _{j}=j\right) \mathbb{P}^{\uparrow }\left(
A_{u.s}|Z_{0}=k\right)  \notag \\
&&\qquad \qquad \qquad =A_{1}(z)\Theta _{j}(K).  \label{Dominate}
\end{eqnarray}

Combining the obtained estimates and letting $K\rightarrow \infty $ we get
\begin{equation*}
\lim_{n\rightarrow \infty }\frac{\mathbb{P}\left( S_{n-m}\leq
za_{m},R_{n}(\varphi (n)),\tau _{n}=j\right) }{\mathbb{P}\left(
Q_{n-j}(\varphi (n))\right) }=A_{1}(z)\Theta _{j}(\infty ),
\end{equation*}%
or, in view of \ (\ref{AsymMain1})
\begin{equation*}
\lim_{n\rightarrow \infty }\frac{\mathbb{P}\left( S_{n-m}\leq
za_{m},S_{n}\leq \varphi (n),Z_{n}>0,\tau _{n}=j\right) }{\mathbb{P}\left(
S_{n}\leq \varphi (n),Z_{n}>0\right) }=A_{1}(z)\frac{\Theta _{j}(\infty )}{%
\Theta }.
\end{equation*}%
Summing over $j$ and taking into account the definition of $\Theta $ in (\ref%
{DefTheta}), \ we get%
\begin{equation*}
\lim_{n\rightarrow \infty }\frac{\mathbb{P}\left( \log Z_{n-m}\leq
za_{m},S_{n}\leq \varphi (n),Z_{n}>0\right) }{\mathbb{P}\left( S_{n}\leq
\varphi (n),Z_{n}>0\right) }=A_{1}(z),
\end{equation*}%
as desired.

2) Assume now that $\varphi (n)\sim Ta_{m},T\in (0,\infty )$. Then,
recalling (\ref{Cond2}) and applying the dominated convergence theorem we
obtain
\begin{eqnarray*}
&&\lim_{n\rightarrow \infty }\frac{\mathbb{P}\left( S_{n-m}\leq
za_{m},R_{n}(\varphi (n)),\tau _{n}=j,S_{j}>-\sqrt{\varphi (n)},Z_{j}\leq
K\right) }{\mathbb{P}\left( Q_{n-j}(\varphi (n))\right) } \\
&&\qquad =\lim_{n\rightarrow \infty }\int_{-\sqrt{\varphi (n)}%
}^{0}\sum_{k=1}^{K}\mathbb{P}\left( S_{j}\in dq,Z_{j}=k,\tau _{j}=j\right) \\
&&\qquad \qquad \times \frac{\mathbb{E}\left[ H_{n-j}(k);S_{n-m-j}\leq
za_{m},Q_{n-j}(\varphi (n)-q)\right] }{\mathbb{P}\left( Q_{n-j}(\varphi
(n))\right) } \\
&&\qquad \qquad =B(z,T)\int_{-\infty }^{0}\sum_{k=1}^{K}\mathbb{P}\left(
S_{j}\in dq,Z_{j}=k,\tau _{j}=j\right) \mathbb{P}^{\uparrow }\left(
A_{u.s}|Z_{0}=k\right) \\
&&\qquad \qquad \qquad =B(z,T)\Theta _{j}(K).
\end{eqnarray*}

Summing over $j$ and taking into account (\ref{AsymMain1}) and the
definition of $\Theta $,\ we get%
\begin{equation*}
\lim_{n\rightarrow \infty }\frac{\mathbb{P}\left( \log Z_{n-m}\leq
za_{m},S_{n}\leq \varphi (n),Z_{n}>0\right) }{\mathbb{P}\left( S_{n}\leq
\varphi (n),Z_{n}>0\right) }=B(z,T),
\end{equation*}%
as desired.

3) Finally, consider the case $m=m(n)\rightarrow \infty $ and $a_{m}=o\left(
\varphi (n)\right) $ as $n\rightarrow \infty $. We introduce the notation $%
S_{n-m,n}:=S_{n-m}-S_{n}$ and for $z\in (-\infty ,\infty )$ write%
\begin{eqnarray}
&&\mathbb{P}\left( \log Z_{n-m}-S_{n}\leq za_{m},R_{n}(\varphi (n))\right)
\notag \\
&=&\mathbb{P}\left( \log \hat{Z}_{n-m}+S_{n-m,n}\leq za_{m},\hat{Z}_{n-m}\in
(M^{-1},M),R_{n}(\varphi (n))\right)  \notag \\
&&+\mathbb{P}\left( \log Z_{n-m}-S_{n}\leq za_{m},\hat{Z}_{n}\notin
(M^{-1},M),R_{n}(\varphi (n))\right) .  \label{Sm1}
\end{eqnarray}%
Now it is not difficult to check that if we replace in all the relations
between (\ref{Decompose0}) and (\ref{LastCommon}) $S_{n-m}$ by $S_{n-m,n}$
and $S_{n-m-j}$ by $S_{n-m-j,n-j}$ then all the estimates between (\ref%
{Decompose0}) and (\ref{LastCommon}) remain valid. As a result,%
\begin{equation*}
\lim_{J\rightarrow \infty }\lim_{n\rightarrow \infty }\frac{\mathbb{P}\left(
\log Z_{n-m}-S_{n}\leq za_{m},R_{n}(\varphi (n)),\tau _{n}>J\right) }{%
\mathbb{P}\left( R_{n}(\varphi (n))\right) }=0
\end{equation*}%
and, for each fixed $j\in \lbrack 0,J]$
\begin{equation*}
\lim_{K\rightarrow \infty }\lim_{n\rightarrow \infty }\frac{\mathbb{P}\left(
\log Z_{n-m}-S_{n}\leq za_{m},R_{n}(\varphi (n)),\tau _{n}=j,Z_{j}>K,S_{j}>-%
\sqrt{\varphi (n)}\right) }{\mathbb{P}\left( R_{n}(\varphi (n))\right) }=0.
\end{equation*}%
Recalling (\ref{Cond3}) and applying the dominated convergence theorem, we
obtain%
\begin{eqnarray*}
&&\lim_{n\rightarrow \infty }\frac{\mathbb{P}\left( S_{n-m,n}\leq
za_{m},R_{n}(\varphi (n)),\tau _{n}=j,S_{j}>-\sqrt{\varphi (n)},Z_{j}\leq
K\right) }{\mathbb{P}\left( Q_{n-j}(\varphi (n))\right) } \\
&&\qquad \qquad \qquad \qquad \qquad \qquad \qquad \qquad \qquad \qquad =%
\mathbf{P}\left( Y_{1}\leq z\right) \Theta _{j}(K).
\end{eqnarray*}%
Combining the estimates above and letting $K$ to infinity, we see that
\begin{eqnarray*}
&&\lim_{n\rightarrow \infty }\frac{\mathbb{P}\left( S_{n-m,n}\leq
za_{m},R_{n}(\varphi (n)),\tau _{n}=j,S_{j}>-\sqrt{\varphi (n)}\right) }{%
\mathbb{P}\left( Q_{n-j}(\varphi (n))\right) } \\
&&\qquad \qquad \qquad \qquad \qquad \qquad \qquad \qquad \qquad \qquad =%
\mathbf{P}\left( Y_{1}\leq z\right) \Theta _{j}(\infty ).
\end{eqnarray*}%
Summing over $j$, taking into account (\ref{AsymMain1}) and the definition
of $\Theta $, we get%
\begin{equation*}
\lim_{n\rightarrow \infty }\frac{\mathbb{P}\left( \log Z_{n-m}-S_{n}\leq
za_{m},S_{n}\leq \varphi (n),Z_{n}>0\right) }{\mathbb{P}\left( S_{n}\leq
\varphi (n),Z_{n}>0\right) }=\mathbf{P}\left( Y_{1}\leq z\right) .
\end{equation*}

Theorem \ref{T_global} is proved.

\textbf{Acknowledgment. }The work of E.E. Dyakonova and  V.A. Vatutin was
performed at the Steklov International Mathematical Center and supported by
the Ministry of Science and Higher Education of the Russian Federation
(agreement no. 075-15-2022-265). The research of C.Dong and V.A. Vatutin was
also supported by the Ministry of Science and Technology of PRC, project
G2022174007L.

\end{document}